\documentclass[12pt,fleqn,a4paper]{amsart}

\usepackage[utf8]{inputenc}
\usepackage[foot]{amsaddr}
\usepackage[pdftex]{lscape}
\usepackage{comment}
\usepackage[in]{fullpage}
\usepackage[onehalfspacing]{setspace}
\usepackage{amsfonts,amssymb,amsthm}
\usepackage[pagebackref,colorlinks,breaklinks,linkcolor=blue,anchorcolor=blue,citecolor=blue]{hyperref}
\usepackage[alphabetic,nobysame]{amsrefs}
\input amsart-modif
\input backref.patch

\numberwithin{equation}{section}\swapnumbers
\usepackage{parskip}

\usepackage[cmtip,matrix,arrow,curve]{xy}
\SelectTips{cm}{10}
\newcommand{\cxymatrix}[1]{\vcenter{\xymatrix@=15pt{#1}}}
\newcommand{\kxymatrix}[1]{\vcenter{\xymatrix@=5pt{#1}}}
\usepackage{rotating}

\newcommand{\xysubseteq}{\ar@{}[r]|{\displaystyle\subseteq}}
\newcommand{\xysubseteqdown}{\ar@{}[d]|{\rotatebox{90}{$\supseteq$}}}

\usepackage{tikz}
\usetikzlibrary{calc}

\usepackage{aliascnt}

\newtheorem{theorem}{Theorem}[section]
\newaliascnt{lemma}{theorem}

\newtheorem{lemma}[lemma]{Lemma}
\aliascntresetthe{lemma}
\newaliascnt{corollary}{theorem}

\newtheorem{corollary}[corollary]{Corollary}
\aliascntresetthe{corollary}
\newaliascnt{proposition}{theorem}

\newtheorem{proposition}[proposition]{Proposition}
\aliascntresetthe{proposition}

\theoremstyle{definition}
\newaliascnt{definition}{theorem}

\newtheorem{definition}[definition]{Definition}
\aliascntresetthe{definition}

\newaliascnt{remark}{theorem}
\newtheorem{remark}[remark]{Remark}
\aliascntresetthe{remark}
\newtheorem{remarks}[remark]{Remarks}

\newtheorem*{remark*}{Remark}

\newaliascnt{example}{theorem}
\newtheorem{example}[example]{Example}
\newtheorem{examples}[example]{Examples}

\aliascntresetthe{example}

\usepackage{enumitem}
\setlist[enumerate,2]{label=\textit{\alph*)},ref=\textit{\alph*}),noitemsep}
\setlist[enumerate,1]{label=\textit{\roman*)},ref=\textit{\roman*}),noitemsep}

\usepackage[capitalize]{cleveref}
\usepackage{microtype}

\renewcommand\[{\begin{equation}}
\renewcommand\]{\end{equation}}

\renewcommand\tilde{\widetilde}

\renewcommand\phi{\varphi}
\renewcommand\epsilon{\varepsilon}

\renewcommand\theta{\vartheta}
\renewcommand\rho{\varrho}

\newcommand\CC{{\mathbb C}}
\newcommand\FF{{\mathbb F}}

\newcommand\KK{{\mathbb K}}
\newcommand\NN{{\mathbb N}}

\newcommand\QQ{{\mathbb Q}}\newcommand\RR{{\mathbb R}}

\newcommand\ZZ{{\mathbb Z}}
\newcommand\cA{{\mathcal A}}
\newcommand\cC{{\mathcal C}}

\newcommand\cR{{\mathcal R}}
\newcommand\cT{{\mathcal T}}

\newcommand\leer{\varnothing}

\newcommand\eq{{\overline{e}}}

\newcommand\Gq{{\overline{G}}}

\newcommand\iq{{\overline{\imath}}}

\newcommand\Kq{{\overline{K}}}

\newcommand\xq{{\overline{x}}}
\newcommand\yq{{\overline{y}}}
\newcommand\zq{{\overline{z}}}
\newcommand\into{\mathrel{\kern-3pt\xymatrix@=10pt{\ar@{>->}[r]&}\kern-5pt}}
\newcommand\inj{\ar@{>->}}
\newcommand\sur{\ar@{>>}}
\newcommand\auf{\twoheadrightarrow}
\newcommand\fau{\twoheadleftarrow}
\newcommand\qur{\mathrel{\kern-10pt\raise2pt\hbox{\xymatrix@=10pt{\ar@{{ }>>}[r]&}\kern-5pt}}}
\newcommand\quot{\mathrel{\kern-5pt\raise2pt\hbox{\xymatrix@=10pt{&\ar@{{ }>>}[l]}\kern-10pt}}}
\newcommand\Auf[1]{\mathrel{\kern-3pt\xymatrix@=10pt{\ar@{>>}[r]^{#1}&}\kern-5pt}}
\newcommand\lrarrow{\mathop{\lower3pt\vbox{\baselineskip0pt\hbox{$\rightarrow$}\hbox{$\leftarrow$}}}}

\renewcommand\*{{\bf1}}
\newcommand\0{{\bf 0}}
\usepackage{ bbold }
\renewcommand\1{{\mathbb{1}\vline width0pt height1pt}}
\newcommand\FS{{\mathsf{Set}}}

\renewcommand\hat{\widehat}
\newcommand\<{\langle}
\renewcommand\>{\rangle}
\def\hide#1{\hbox to 10pt{\hss$#1$\hss}}
\def\|#1|{\operatorname{#1}}
\newcommand\sO{{\mathrm s}}
\newcommand\qO{{\mathrm q}}
\newcommand\sqO{{\mathrm{sq}}}
\newcommand{\doppelpfeil}{\rightrightarrows}
\makeatletter
\newcommand{\Times}{\gen@tens{\times}}
\newcommand{\gen@tens}[1]{%
  \@ifnextchar_{\gen@@tens{#1}}{\mathbin{#1}}%
}
\def\gen@@tens#1_#2{%
  \mathpalette\gen@@@tens{{#1}{#2}}%
}
\newcommand\gen@@@tens[2]{\mathbin{\gen@@@@tens#1#2}}
\newcommand\gen@@@@tens[3]{%
  \ifx#1\displaystyle
    \mathop{#2}\limits_{#3}%
  \else
    {#2}_{#3}%
  \fi
}
\makeatother
\newcommand{\sfrac}{\textstyle\frac}
\newcommand\ueber[2]{\genfrac{}{}{0pt}{}{#1}{#2}}

\newcommand{\Hide}[1]{\hbox to 0pt{\hss$#1$\hss}}
\newcommand{\St}{h}
\newcommand\x{{\underline x}}
\newcommand\y{{\underline y}}
\newcommand\z{{\underline z}}

\newcommand\9{''}

\newif\ifcharacterstart
\characterstarttrue

\title[Multiplicities and dimensions]{Multiplicities and dimensions in enveloping\\
  tensor categories}
\author[]{Friedrich Knop$^1$}
\address[]{$^1$Department Mathematik\\FAU Erlangen-Nürnberg\\
  Cauerstraße 11\\
  D-91058 Erlangen\\
  Germany}
\email{friedrich.knop@fau.de}
\date{November 27, 2023}

\subjclass[2020]{Primary 18M05, 18B10, 18E13; Secondary 05E05, 18E08, 20C15}

\keywords{Tensor categories; Semisimple categories; Regular
  categories; Mal’cev categories; Möbius functions; Lattices; Symmetric functions; Finite groups}

\begin{document}

\begin{abstract}
  In a previous paper, semisimple tensor categories were constructed
  from certain regular Mal'cev categories. In this paper, we
  calculate the tensor product multiplicities and the categorical
  dimensions of the simple objects. This yields also the Grothendieck
  ring. The main tool is the subquotient decomposition of the generating
  objects.
\end{abstract}

\maketitle         

\section{Introduction}

Fix an algebraically closed field $\KK$ of characteristic $0$. By a
\emph{tensor category} we mean a $\KK$-linear, locally finite, abelian
category which is equipped with a symmetric tensor product such that
each object has a dual object (rigidity) and such that
$\|End|\1=\KK$. They are also referred to as \emph{symmetric tensor
  categories} or \emph{pre-Tannakian categories}.  A typical example
is the category of finite dimensional $\KK$-vector spaces or, more
general, representations of a linear algebraic group. Categories of
this type are called Tannakian.

First examples of non-Tannakian categories were found by Deligne, see
e.g., \cite{Deligne}, using an interpolation procedure.  In
\cite{TERC} a different method for the construction of non-Tannakian
categories was devised. It is based on the observation that the
category of relations of a base category $\cA$ is already rigid,
symmetric, and monoidal. After making it $\KK$-linear and twisting the
product of relations by a $\KK$-valued degree function $\delta$ one
obtains a category $\cT=\cT(\cA,\delta)$ which is in many cases
semisimple, hence in particular abelian.

The necessary assumptions on $\cA$ and $\delta$ recalled in
\cref{sec:RCTE}. They hold, for example, for any category of finite
algebraic objects containing a group structure, like finite groups,
finite rings, finite modules, Boolean algebras and many more.
Deligne's category from \cite{Deligne} is recovered by taking for
$\cA$ the opposite category of the category of finite sets.

In the present paper we start analyzing the internal structure of
$\cT$ in case it is semisimple. Besides $\cT$ being an interesting
object in its own right there is another motivation: Despite $\cT$ not
being an interpolation category \emph{a priori}, in many cases it can
be realized as one \emph{a posteriori}. Therefore, it may be possible
to transfer properties from $\cT$ to the interpolated categories. This
paper contains two examples for this transfer.

As for any semisimple category the most important objects of $\cT$ are
the simple ones. In fact, these were already determined in \cite{TERC}
as part of the semisimplicity proof: simple objects of $\cT$ are
classified by pairs $(x,\chi)$ where $x$ is an object of $\cA$ (up to
isomorphism) and $\chi$ is an irreducible character of its
automorphism group. The corresponding simple object is denoted by
$[x]^0_\chi$.

It is now a natural problem how the simple objects relate to the
tensor product. There are two classical questions: \emph{a)} how does
the tensor product $[x]^0_\chi\otimes [y]^0_\eta$ decompose in simple
objects and \emph{b)} what is the (internal) dimension
$\dim_\cT[x]^0_\chi$ of a simple object? Both of these questions are
being answered in the present paper.

For \emph{i)} define for any three objects $x_1$, $x_2$, and $x_3$ the
set
  \[
    T(x_1,x_2,x_3):=\{r\subseteq x_1\times x_2\times x_3\mid
    r\auf x_1,x_2,x_3; r\into x_1\times x_2,x_1\times x_3,x_2\times x_3\}.
  \]
  Let $\chi_{T(x_1,x_2,x_3)}$ be the induced permutation character of
  $\|Aut|_\cT(x_1)\times\|Aut|_\cT(x_2)\times\|Aut|_\cT(x_3)$.

  \begin{theorem}[\cref{cor:multtensor2} below]
    Let $[x_1]^0_{\chi_1}$ and $[x_2]^0_{\chi_2}$ be two simple
    objects of $\cT$. Then
  \[
    [x_1]^0_{\chi_1}\otimes[x_2]^0_{\chi_2}\cong
    \bigoplus_{x,\chi}\ \<\chi_T(x_1,x_2,x)\mid
    \chi_1\otimes\chi_2\otimes\chi^\vee\>\ [x]^0_{\chi}.
   \]
\end{theorem}

For the dimension formula we give two variants. For the simpler one we
need the following notation. For any object $x$ let $\sO(x)$ and
$\qO(x)$ the lattice of subobjects and quotient objects, respectively,
where $\qO(x)$ is ordered in such a way that $z=x$ is the minimum. For
any automorphism $g$ of $x$ let $\qO(x)^g=\{z\in\qO(x)\mid
gz=z\}$. Let $x_g\in\qO(x)^g$ be the minimal element $z$ such that
$g|_z=\|id|_z$. Let moreover $\hat x_g\in\qO(x)^g$ be the join of all
atoms. Finally, let $\mu_\Lambda$ denote the Möbius function of a
lattice $\Lambda$. Define, in particular,
\[
  \omega_x:=\sum_{y\in\sO(x)}\mu_{\sO(x)}(y,z)\ \delta(y\auf\*).
\]

\begin{theorem}[\cref{cor:second} below]
  Let $x$ be an object of $\cA$ and $\chi$ an irreducible character of
  $A=\|Aut|_\cT(x)$. Then
  $\dim_\cT[x]_\chi^0=\<\,\chi_{[x]^0}\mid\chi\,\>_A$ where the class
  function $\chi_{[x]^0}$ on $A$ is defined as
  \[
    \chi_{[x]^0}(g):=\sum_{x_g\le z\le\hat x_g}
    \mu_{\qO(x)^g}(x,z)\ \omega_z
  \]
\end{theorem}

These simple formulas become quite complicated when interpreted in a
concrete category $\cA$. We succeeded in doing so for the case
$\cA=\mathsf{Set}^{\|op|}$. This way, we were able to rederive a
formula of Littlewood \cite{Littlewood} for stable Kronecker
coefficients (the equivalence is proved in Appendix B) and to find an
apparently new identity for symmetric functions (for which we found a
symmetric function proof a posteriori, see also Appendix B).

The category $\cT$ is built from basic objects $[x]$ where $x$ is an
object of $\cA$. So the main technical result of this paper is the
\emph{subquotient decomposition} which also yields a decomposition of
$[x]$ into simple objects. This decomposition builds up from the
coarser \emph{subobject decomposition} of our previous paper
\cite{SOD}.

The methods of the paper are mostly category theoretical but we also
use a fair amount of lattice theory and symmetric function theory. To
not interrupt the flow of the paper we collected the necessary facts
in two appendices.

\section{Regular categories and their tensor envelopes}
\label{sec:RCTE}

In this section we briefly recall terminology, notation, and the
construction of $\cT(\cA,\delta)$. For details see \cites{TERC,SOD}.

We start with a category $\cA$ and use the following terminology:
Monomorphisms will be called \emph{injective} and denoted by
$y\into x$. Every injective morphism defines a \emph{subobject} of its
target. The \emph{image} of a morphism $y\to x$ is the smallest
subobject of $x$, the morphism factors through. Morphisms whose image
equals the target (i.e., extremal epimorphisms) will be called
\emph{surjective} and denoted by $x\auf z$.

We assume throughout that $\cA$ is \emph{finitely
  complete}\footnote{This condition can be slightly relaxed to
  accommodate, e.g., the category of affine spaces \cite{TERC}.}. The
terminal object will be denoted by $\*$.  We also assume that $\cA$ is
\emph{regular}. This means that all morphisms have images and that
surjectivity is preserved under pullback.

A \emph{relation} between objects $x$ and $y$ of $\cA$ is a subobject
$r$ of $x\times y$. If $s\subseteq y\times z$ is another relation then
their \emph{product} $r\circ s$ is the relation
$r\circ s=\|image|(r\times_ys\to x\times z)$. In regular categories,
the product of relations is associative. This way, we obtain the
\emph{category $\|Rel|(\cA)$ of relations} of $\cA$ with the same
objects but relations as morphisms.

The category $\cT(\cA,\delta)$ is a twisted $\KK$-linear version of
$\|Rel|(\cA)$. The twist function $\delta$ is a map which assigns to
any epimorphism $e$ an element $\delta(e)$ of some fixed base field
$\KK$. It is subject to the requirements {\bf D1}--{\bf D3} from
\cite{TERC}.

  
The category $\cT^0=\cT^0(\cA,\delta)$ will have the same objects as
$\cA$. More precisely, each object $x$ of $\cA$ gives rise to an
object $[x]$ of $\cT$. The morphisms $[x]\to[y]$ are the formal
$\KK$-linear combinations of relations $r\subseteq x\times y$. The
morphism induced by the relation $r$ will be denoted by $\<r\>$. The
product of two such morphism is the $\delta$-twisted product of
relations:
\[
  \<r\>\<s\>:=\delta(r\times_ys\auf r\circ s)\ \<r\circ s\>.
\]

The category $\cT(\cA,\delta)$ is the \emph{pseudo-abelian completion}
of the \emph{additive completion} of $\cT^0$, i.e., the category
obtained by formally adjoining direct sums and direct summands.

With the tensor product induced by the direct product,
$[x]\otimes[y]=[x\times y]$, $\cT(\cA,\delta)$ becomes a
pseudo-abelian, symmetric, monoidal category. Its unital object is
$\1=[\*]$. It is also rigid with $[x]^\vee=[x]$ for all $x$. The
adjoint of $f=\<r\>$ is $f^\vee=\<r^\vee\>$ where $r^\vee$ is the
relation obtained by swapping the factors. Every morphism $f:x\to y$
gives rise to a morphism $[f]=\<\Gamma_f\>:[x]\to[y]$ where
$\Gamma_f=\|im|(id\times f:x\to x\times y)$ is the graph of $f$. This
way, one obtains a faithful functor $\cA\to\cT$.
  
For $\cT$ to be semisimple more assumptions on $\cA$ are needed.
First, $\cA$ has to be \emph{subobject finite}, i.e., every object has
only finitely many subobjects. This condition ensures that all
morphism spaces $\|Hom|_\cT(x,y)$ are finitely generated
$\KK$-modules.

Another assumption we are going to make is that $\cA$ is
\emph{exact}. Here exactness (in the sense of Barr \cite{Barr}) means
that all equivalence relations $r\subseteq x\times x$ are effective,
i.e., there always exists a quotient object $x/r$.
  
The most restrictive condition is the requirement of $\cA$ being
Mal'cev. This means, e.g., that every reflexive relation
$r\subseteq x\times x$ is an equivalence relation. The Mal'cev
property enters mostly through the validity of via \emph{Goursat's
  lemma} in $\cA$ (see \cite{TERC}*{(5.1)}. It says:

\emph{In an exact Mal'cev category every relation
  $r\subseteq x\times y$ is of the form $\xq\times_c\yq$ where
  $\xq\subseteq x$ and $\yq\subseteq y$ are subobjects and $c$ is a
  quotient object of both $\xq$ and $\yq$:
  \[\label{eq:dreieck3}
    \cxymatrix{&&r\sur[dl]\sur[dr]\\
      &\xq\inj[dl]\sur[dr]&&\yq\inj[dr]\sur[dl]\\
      x&&c&&y}
  \]
} In short, this means that subobjects of $x\times y$ can be described
in terms of $x$ and $y$ alone.

Despite the restrictiveness of the Mal'cev condition there is an
abundance of important examples. In particular, all algebraic theories
containing a group operation are exact Mal'cev. This includes the
categories of groups, vector spaces, rings, modules, etc. All abelian
categories are exact Mal'cev, as well. Combined with our condition of
subobject finiteness, we are dealing with the finite models of such
theories, like finite groups, and so on.

There is also one example which is slightly less obvious, namely the
category $\mathsf{Set}^{\|op|}$ which is opposite to the category of
finite sets. This is the category which gives rise to Deligne's
category $\|Rep|S_t$ of \cite{Deligne}. It is easily recognized as
exact Mal'cev, though, by observing that it is equivalent to the
category of finite Boolean algebras.

There are two more conditions. First, we will assume throughout $\KK$
is algebraically closed of characteristic zero. Secondly, it is
convenient to assume that the terminal object $\*$ has no proper
subobjects. This is equivalent to $\|End|_\cT(\*)=\KK$. In this case
we will extend the definition of $\delta$ to objects by
\[
  \delta(x):=\delta(x\auf\*).
\]
The condition is quite innocuous since the general case can be reduced
to it (see \cite{TERC}*{Thm.~3.6}).

Now we recall a numerical criterion for $\cT(\cA,\delta)$ to be a
semisimple tensor category. Since $\cA$ is subobject finite, its
collection $\sO(x)$ of subobjects forms a finite set. It is partially
ordered by inclusion with $x$ being the maximum. In particular, the
Möbius number $\mu_\sO(y,x)$ is defined for every $y\in\sO(x)$. To
every surjective morphism $e:x\auf z$ one assigns the following
element of $\KK$:
\[\label{eq:omegax}
  \omega_e:=\sum_{\ueber{y\subseteq x}{e(y)=z}}\mu_\sO(y,x)\delta(e|_y).
\]
When $\sO(\*)=\{\*\}$ we extend $\omega_e$ to objects by setting
\[\omega_x:=\omega_{x\auf\*}.
\]

It is a nontrivial fact, \cite{TERC}*{Lemma~8.4}, that $\omega_e$ is
multiplicative in $e$, i.e., for all
$x\overset e\auf y\overset f\auf z$ one has
\[
  \omega_{fe}=\omega_f\omega_e.
\]
This reduces the computation of $\omega_e$ to the case when $e$ is
indecomposable.

We call $\delta$ \emph{non-degenerate} if $\omega_e\ne0$ for all
(indecomposable) surjective morphisms $e$. Then the criterion is:

\begin{theorem}[\cite{TERC}*{Thm.~6.1 and Thm.~8.3}]\label{thm:main1}
  Let $\cA$ be a regular, subobject finite, exact, Mal'cev category
  and let $\delta$ be a degree function on $\cA$ with values in a
  field $\KK$ of characteristic $0$. Then $\cT(\cA,\delta)$ is a
  semisimple tensor category if and only if $\delta$ is
  non-degenerate.
\end{theorem}

\begin{examples}
  There are three examples which one might want to
  keep in mind when reading the paper.

  \emph{i)} The category
  $\cA=\mathsf{Set}^{\|op|}\cong\mathsf{Bool}$. Here the surjective or
  injective morphisms are the injective or surjective maps between
  finite sets, respectively. The general degree function applied to an
  injective map $e:A\into B$ is $\delta(e)=t^{|B|-|A|}$ where
  $t\in\KK$ is any element. The morphism $e$ is indecomposable if
  $|B|=|A|+1$. In that case $\omega_e=t-|A|$. Hence $\cT(\cA,\delta)$
  is semisimple if and only if $t\not\in\NN$.
  
    \emph{ii)} The category $\cA=\mathsf{Vect}(\FF_q)$. In this case,
    the degree function is $\delta(e)=t^{\dim V-\dim U}$ where
    $e:U\auf V$ is surjective. The morphism $e$ is indecomposable if
    $\dim V=\dim U+1$ in which case $\omega_e=t-q^{\dim U}$. It
    follows that $\cT(\cA,\delta)$ is semisimple if and only if $t$ is
    not of the form $q^n$ with $n\in\NN$. Observe that $t=0$ is
    permitted.
  
    \emph{iii)} The category $\cA=\mathsf{Group}$ of finite groups. In
    this case, the general degree function has infinitely many free
    parameters $t_S$, one for each simple group $S$. The degree
    function then takes up the form $\delta(G)=t_{C_1}\ldots t_{C_n}$
    where $C_1,\ldots,C_n$ are the composition factors of $G$ counted
    with multiplicity. A surjective homomorphism $e:G\auf\Gq$ is
    indecomposable if its kernel $N$ is a minimal, non-trivial, normal
    subgroup. The computation of $\omega_e$ involves all subgroups
    $H\subseteq G$ with $HN=G$. This is doable in every given case but
    impossible in general even if $G=N$ is simple and $\Gq=1$. Since
    $H=G$ is always present one can say at least that
    $\omega_e=\delta(G)+\text{lower order terms}$. In any case,
    $\cT(\cA,\delta)$ will be semisimple if all variables $t_S$ are
    algebraically independent over $\QQ$.
\end{examples}

Let $S$ be an object of $\cT$ whose endomorphism ring is a division
algebra. Then the non-degeneracy of $\delta$ ensures that every
nonzero $\cT$-morphism $X\to S$ onto $S$ has a section. This and
induction on the size of $x$ is used to prove the following lemma:

\begin{lemma}\label{lemma:defx1}
  Under the assumptions of \cref{thm:main1} let $x$ be an object of
  $\cA$. Then
  
  \begin{enumerate}

  \item\label{it:defx1-1} $[x]$ has a unique direct summand $[x]^1$
    with:

    \begin{enumerate}
    \item\label{it:defx1-1-1} $[x]^1$ is isomorphic to a direct
      summand of an object of the form $[z_1]\oplus\ldots\oplus[z_n]$
      where each $z_i$ is proper subquotient of $x$.

    \item\label{it:defx1-1-2} Conversely, every morphism $[z]\to[x]$
      factors through $[x]^1$ whenever $z$ is a proper subquotient of
      $x$

    \end{enumerate}

  \item\label{it:defx1-2} $[x]^1$ has a unique complement $[x]^0$ in
    $[x]$.

  \item\label{it:defx1-3} $\KK[\|Aut|_\cA(x)]\to\|End|_\cT([x]^0)$ is
    an isomorphism.

  \end{enumerate}

\end{lemma}

The proof of part \emph{iii)} makes crucial use is made of Goursat's
lemma. It implies that any relation $r\subseteq x\times x$ is either
the graph of an automorphism or the induced morphism $\<r\>:[x]\to[x]$
factors through $[c]$ where $c$ is a proper subquotient of $x$ (see
diagram \eqref{eq:dreieck3} with $y=x$).

Part \emph{iii)} also yields a classification of the simple objects of
$\cT$. Let $A:=\|Aut|_\cA(x)$. Then $[x]^0$ has an $A$-isotypical
decomposition
\[\label{eq:decomp}
  [x]^0=\bigoplus_\chi V_\chi\otimes_\KK [x]^0_\chi
\]
where $V_\chi$ is the irreducible $A$-module with character $\chi$ and
$[x]^0_\chi$ is a simple object of $\cT$. Conversely, every simple
object of $\cT$ is uniquely of that form.

\section{The subobject decomposition}

In this section $\cA$ may be any regular, subobject finite category
with arbitrary degree function $\delta$.  Our goal is to decompose any
object $[x]$ as far as possible as a direct sum. We briefly recall the
subobject decomposition which was studied in detail in \cite{SOD}.

Any decomposition requires the construction of idempotents. As a
starter, each subobject $y$ of $x$ gives rise to the relation which is
just the diagonal embedding $y\into x\times x$. This defines an
endomorphism $p_y=\<y\>$ of $[x]$ which is easily seen to be
idempotent.

Now, the $p_y$ for different $y$ commute with each other. Hence all
decompositions $[x]=p_y[x]+(1-p_y)[x]$ have a common refinement
\[\label{eq:SOD}
  [x]=\bigoplus_{y\subseteq x}[y]^*,
\]
called the \emph{subobject decomposition} of $[x]$. Here
$[y]^*:=p_y^*[x]$ where $p_y^*$ is a primitive idempotent in the
algebra generated by the $p_y$. These are computed by Möbius
inversion:
\[\label{eq:p*}
  p_y^*:=\sum_{y'\subseteq y}\mu_\sO(y',y)p_{y'}\quad\text{and}\quad
  p_y=\sum_{y'\subseteq y}p_{y'}^*.
\]
To be very precise, the subobjects $[y]^*\subseteq[x]$ also depend on
$x$ but it is easy to show that the relation $y\into y\times x$
defines an isomorphism between $[y]^*\subseteq[y]$ and
$[y]^*\subseteq[x]$. So we are not distinguishing between these two.

It follows from \eqref{eq:SOD} that one can take the objects $[x]^*$
also as generators of the pseudo-abelian tensor category $\cT$.  In
fact, Deligne's construction of
$\|Rep|S_t=\cT(\mathsf{Set}^{\|op|},\delta)$ does exactly that. The
advantage of this route is that the space $\|Hom|_\cT([x]^*,[y]^*)$ is
much smaller than $\|Hom|_\cT([x],[y])$. On the other hand, the
definitions of compositions, tensor products, and even the
associativity constraint become much more involved. This approach has
been carried out in the paper \cite{SOD}.

\textbf{For the remainder of this paper we make the following
  assumptions on $\cA,\boldsymbol\delta$:}

\emph{1. The assumptions of \cref{thm:main1} hold, i.e., $\cA$ is a
  regular, subobject finite, exact, Mal'cev category and $\delta$ is a
  non-degenerate degree function on $\cA$ with values in a field $\KK$
  of characteristic $0$. This implies that $\cT=\cT(\cA,\delta)$ is a
  semisimple tensor category.}

\emph{2. Additionally, we assume for simplicity that $\*$ has no
  proper subobject or, equivalently, that $\|End|_\cT(\1)=\KK$ and
  that $\KK$ is algebraically closed.}

\section{The subquotient decomposition}

We have seen that every subobject $y\in\sO(x)$ gives rise to an
idempotent $p_y$. What about quotients $e:x\auf z$? The obvious
candidate would be the relation $r=x\times_zx\subseteq x\times x$. It
defines an endomorphism $q_z=\<r\>=[e]^\vee[e]$ of $[x]$ which is
almost idempotent in that $q_z^2=\delta(e)q_z$. Hence
$\delta(e)^{-1}q_z$ is an idempotent on $[x]$ provided
$\delta(e)\ne0$. Unfortunately, $\delta(e)$ may vanish even if
$\delta$ is non-degenerate. A more serious disadvantage of the
idempotents $q_z$ is the fact that they do not commute with the
idempotents $p_y$.

But still, one can embed $[z]$ into $[x]$ in a canonical way. It is
just not defined by a canonical idempotent.

\begin{lemma}\label{lem:surj}
  Let $e:x\auf z$ be surjective. Then $[e]^\vee:[z]\to[x]$ is
  injective.
\end{lemma}

\begin{proof}
  It suffices to show that the map
  $[e]^\vee\circ:\|End|_\cT([z])\to\|Hom|_\cT([z],[x])$ is
  injective. Indeed, since $[z]$ is semisimple, any projection onto
  the kernel of $[e]^\vee$ would be in the kernel of $[e]^\vee\circ$.

  Let $r\into z\times z$ be a relation. Then $\<r\>\in\|End|_\cT([z])$
  is a basis element and we have
  \[
    [e]^\vee\circ\<r\>=
    \left\<\cxymatrix{&&\hide{r\times_zx}\ar@{..>>}[ld]_{e'}\ar@{..>}[rd]^{b'}\\
        &r\ar[ld]_a\ar[rd]^b&&x\sur[ld]_e\ar@{=}[rd]&\\z&&z&&x}\right\>
    =
    \left\<\cxymatrix{&\hide{r\times_zx}\ar[ld]_>>>{ae'}\ar[rd]^>>>{b'}\\
        z&&x}\right\>.
  \]
  Indeed, the right equality holds because the morphism
  $r\times_zx\to z\times x$ is injective which can be seen with the
  commutative diagram
  \[
    \xymatrix{r\times_zx\ar[r]^{(ae',b')}\inj[d]_{(e',b')}&z\times x\ar[d]^{(\|id|_z,e\times\|id|_x)}\\
      r\times x\inj[r]_{(a\times b,\|id|_x)}&z\times z\times x}
  \]
  Since $r$ can be recovered from $r\times_zx$ as image of $(ae',eb')$
  in $z\times z$ the homomorphism $[e]^\vee\circ$ induces an injective
  map from a basis into a basis. Hence it is injective itself.
\end{proof}

Next we combine subobjects and quotient objects. Recall that $z$ is a
\emph{subquotient} of $x$ if there exists a subobject $y\into x$ and a
surjection $y\auf z$. When we remember $y$ we talking about a
\emph{subquotient object}. More precisely, a subquotient object of $x$
is an equivalence class of diagrams $z\fau y\into x$ where two
diagrams are equivalent if there are (necessarily unique) isomorphisms
rendering the following diagram commutative:
\[
  \cxymatrix{
    z\ar[d]^\sim&y\sur[l]\inj[r]\ar[d]^\sim&x\ar@{=}[d]\\
    z'&y'\sur[l]\inj[r]&x }
\]
We denote this subquotient object also by $z_y$. The set of
subquotient objects of $x$ will be denoted by $\sqO(x)$. It is easy to
see that any subquotient $u$ of a subquotient $z$ of $x$ is also a
subquotient of $x$ (see, e.g., \eqref{eq:subsub} below).

A \emph{quotient object} is a subquotient object $z_y$ with $y=x$.
That is the notion dual to subobject. Let $\qO(x)$ be the set of
quotient objects of $x$. Like $\sO(x)$ also $\qO(x)$ is partially
ordered with $z\le z'$ if $x\auf z\auf z'$. Thus $z=x$ is the minimum
of $\qO(x)$.

Any subquotient object $z\overset e\fau y\overset i\into x$ induces
the morphism
\[
  \iota:[i][e]^\vee:[z]\into[x].
\]
which is injective by \cite{SOD}*{(4.1)} and \cref{lem:surj}. Denote
its image by $[z]_y$. It depends only on the subquotient object
represented by $e,i$. Moreover, $\iota:[z]\to[z]_y$ is an isomorphism.

\begin{lemma}\label{lem:x1-neu}
  For an object $[x]$ of $\cA$ let $[x]^1$ as in
  \cref{lemma:defx1}. Then
  \[
    [x]^1=\sum_{\ueber{z_y\in\sqO(x)}{ z_y\ne x_x}}[z]_y.
  \]
\end{lemma}
  
\begin{proof}
  We show that the right hand side $[x]'$ satisfies the
  characterization of $[x]^1$ from
  \cref{lemma:defx1}\ref{it:defx1-1}. There is an surjective morphism
  $\bigoplus_{\ueber{z_y\in\sqO(x)}{ z_y\ne x_x}}[z]\auf[x]'$. Since $\cT$
  is semisimple it has a section. Hence part \ref{it:defx1-1-1} is
  satisfied. For part \ref{it:defx1-1-2} let $f:[z]\to[x]$ be a
  morphism with $z$ being a proper subquotient of $x$. We may assume
  that $f=\<r\>$ for some relation $r\subseteq z\times x$. By
  Goursat's lemma we get a diagram
  \[\label{eq:dreieck}
    \cxymatrix{&&r\sur[dl]\sur[dr]\\
      &\zq\inj[dl]\sur[dr]&&\xq\inj[dr]\sur[dl]\\
      z&&c&&x}
  \]
  Since $c$ is a subquotient of $z$ it is a proper subquotient of
  $x$. Hence $f([z])\subseteq [c]_\xq\subseteq[x]'$.
\end{proof}

Recall that $[x]^1$ has a unique complement $[x]^0$ in $[x]$. Then
each $z_y\in\sqO(x)$ induces the subobject
\[
  [z]^0_y:=\iota[z]^0\subseteq[x].
\]

\begin{theorem}[Subquotient decomposition]\label{thm:subquotient}
  Let $[x]$ be an object of $\cA$. Then
  \[\label{eq:subquotient1}
    [x]=\bigoplus_{z_y\in\sqO(x)}[z]_y^0.
  \]
  This decomposition is compatible with the subobject
  \emph{filtration}, i.e.,
  \[\label{eq:subquotient2}
    p_y[x]=\bigoplus_{\ueber{z_u\in\sqO(x)}{ u\subseteq y}}[z]_u^0.
  \]
  Moreover, the projection $p_x^*:[x]\to[x]^*$ induces an isomorphism
  \[\label{eq:subquotient3}
    \bigoplus_{z\in\qO(x)}[z]_x^0\overset\sim\to[x]^*.
  \]
\end{theorem}

\begin{proof}
  Let $X:=\bigoplus_{z_y\in\sqO(x)}[z]_y^0$ be the right hand side of
  \eqref{eq:subquotient1}. We claim that the canonical morphism
  $\Phi_x:X\to[x]$ is surjective. To this end we use induction on
  $|\sqO(x)|$. Thus we assume that $\Phi_z$ is surjective for all
  proper subquotients $z$ of $x$. Let $z_y$ be a subquotient object of
  $x$ and $\zq_\yq$ a subquotient object of $[z]$. Then we obtain the
  following diagram with Cartesian square:
  \[\label{eq:subsub}
    \cxymatrix{y'\ar@{>.>}[r]^{i'}\ar@{.>>}[d]_{e'}&y\inj[r]^i\sur[d]_e&x\\
      \yq\inj[r]^\iq\sur[d]_\eq&z\\
      \zq}
  \]
  Thus $\zq_{y'}\in\sqO(x)$. Let $\iota_\yq:[\zq]\to[z]$,
  $\iota_y:[z]\to[x]$, and $\iota_{y'}:[\zq]\to[x]$ the induced
  morphisms. Because of $[e]^\vee[\iq]=[i'][e']^\vee$ (easy, see
  \cite{SOD}*{\bf Rel2}) we have
  $\iota_y\circ\iota_\yq=\iota_{y'}$. Hence \cref{lem:x1-neu} implies
  \[
    \begin{split}
      [x]=[x]^0+[x]^1&=[x]^0+\sum_{\ueber{z_y\in\sqO(x)}{
                       z_y\ne x_x}}\iota_y[z]=\\
                     &=[x]_x^0+\sum_{\ueber{z_y\in\sqO(x)}{
                       z_y\ne x_x}}\iota_y\sum_{\zq_\yq\in\sqO(z)}\iota_\yq[\zq]^0
                       =\sum_{\zq_{y'}\in\sqO(x)}\iota_{y'}[\zq]^0=\Phi_x(X),
    \end{split}
  \]
  proving our claim.
  
  Now we claim that
  $\dim_\KK\|End|_\cT([x])=\dim_\KK\|End|_\cT(X)$. Because $X$ is
  semisimple this will prove that $\Phi_x$ is an isomorphism.  The
  dimension of $\|End|_\cT([x])$ is, by definition, the number of
  subobjects $r$ of $x\times x$. By Goursat's lemma, each $r$ can be
  represented as a diagram
  \[
    \cxymatrix{&&r\sur[ld]\sur[rd]\\
      &\xq\inj[ld]\sur[rd]&&\xq'\sur[ld]\inj[rd]&\\
      x&&z&&x}
  \]
  which is unique up to isomorphism. Thus subobjects of $x\times x$
  are classified by triples $(s_1,s_2,\phi)$ where
  $s_1=z_\xq,s_2=z_{\xq'}\in\sqO(x)$ are subquotient objects of $x$
  and $\phi:z\to z$ is an isomorphism between the corresponding
  subquotients. Let $\hat\cA$ denote a set of representatives of the
  isomorphism classes of objects of $\cA$ and, for $z_0\in\hat\cA$,
  let $\sqO_{z_0}(x)$ be the set of subquotient objects $z_y$ with
  $z\cong z_0$. Then
  \[
    \dim_\KK\|End|_\cT[x]=|\sO(x\times
    x)|=\sum_{z\in\hat\cA}|\sqO_z(x)|^2\ |\|Aut|_\cA(z)|.
  \]
  On the other side we have $\|Hom|([z_1]^0,[z_2]^0)=0$ unless
  $z_1\cong z_2$. Thus,
  \[
    \dim_\KK\|End|_\cT(X)=\sum_{z\in\hat\cA}|\sqO_z(x)|^2\
    \dim_\KK\|End|_\cT([z]^0)=\sum_{z\in\hat\cA}|\sqO_z(x)|^2\
    |\|Aut|_\cA(z)|.
  \]
  This proves the claim and therefore \eqref{eq:subquotient1}. Then
  \eqref{eq:subquotient2} follows from \eqref{eq:subquotient1} with
  $x$ replaced by $y$ and the fact that
  $\{z_u\in\sqO(x)\mid u\subseteq y\}=\sqO(y)$. The isomorphism
  \eqref{eq:subquotient3} follows immediately by observing that
  $[x]^*\cong[x]/\sum_{y\subsetneq x}[y]$.
\end{proof}

\begin{remark}
  The decomposition \eqref{eq:subquotient3} combined with the
  subobject decomposition \eqref{eq:SOD} yields
  \[
    [x]=\bigoplus_{y\subseteq x}[y]^*\cong\bigoplus_{z_y\in\sqO(x)}[z]^0.
  \]
  This decomposition is \emph{not} the same as the subquotient
  decomposition \eqref{eq:subquotient1}. It is rather its associated
  graded. This is because in general $[z]_x^0\not\subseteq[x]^*$.
\end{remark}

\section{Tensor product multiplicities}\label{section:tensor product}

In this section we derive a formula for the tensor product of any two
simple objects of $\cT$ and thereby compute the Grothendieck ring. In
this section we use the following notation: Given $n\ge0$ objects
$x_1,\ldots,x_n$ of $\cA$ we put
\[
  \x:=\prod_j x_j\text{ and }x^i:=\prod_{j\ne i}x_i.
\]
Moreover, let
\[
  T(x_1,\ldots,x_n):=\{r\subseteq\x\mid r\auf x_i\text{ and }r\into
  x^i\text{ for all }i\}.
\]

\begin{remark}
  For $n=0$ we have $\x=\*$ and therefore $T(\emptyset)=\{\*\}$. For
  $n=1$ we have $x^1=\*$ and therefore $T(\*)=\{\*\}$ and
  $T(x_1)=\leer$ when $x_1\not\cong\*$. For $n=2$ we have
  $x^i=x_{3-i}$. Hence for $r$ to lie in $T(x_1,x_2)$ all projections
  $r\to x_i$ must be isomorphisms. Thus $T(x_1,x_2)$ consists of
  graphs of isomorphisms $x_1\overset\sim\to x_2$.
\end{remark}

\begin{theorem}
  Let $x_1,\ldots,x_n$ be objects of $\cA$. Then there is a canonical
  isomorphism
  \[\label{eq:abc2}
    \phi:\<T(x_1,\ldots,x_n)\>_\KK\overset\sim\to
    \|Hom|_\cT(\1,[x_1]^0\otimes\ldots\otimes[x_n]^0).
  \]
\end{theorem}

\begin{proof}
  By definition there is an isomorphism
  \[\label{eq:phiiso}
    \phi:\<\sO(x)\>_\KK\overset\sim\to \|Hom|_\cT(\1,[x]).
  \]
  Using $[x_i]=[x]^0\oplus\sum_{z_y\ne (x_i)_{x_i}}[z]_y$ (see
  \cref{lemma:defx1} and \cref{lem:x1-neu}) one gets
  \[
    [x]=\bigotimes_i[x_i]^0\oplus X
  \]
  where $X$ is the sum of all $[\z]_\y\subseteq[x]$ of the form
  $\z=\prod_iz_i$, $\y=\prod_iy_i$ such that at least one subquotient
  object $(z_i)_{y_i}$ of $x_i$ proper.

  Let $T'$ be the complement of $T(x_1,\ldots,x_n)$ in $\sO(\x)$. We
  claim that $\phi(T')\subseteq\|Hom|_\cT(\1,X)$. Fix $r\in T'$ and
  let $y_i\subseteq x_i$ and $r_i\subseteq x^i$ be the images of
  $r$. Goursat's lemma applied to $r\subseteq\x=x^i\times x_i$ yields
  a diagram
  \[\label{eq:xizixi}
    \cxymatrix{
      &&r\sur[dl]_a\sur[dr]\\
      &r_i\inj[dl]\sur[dr]&&y_i\sur[dl]_{a'}\inj[dr]^b\\
      x^i&&z_i&&x_i }
  \]
  where the square is Cartesian. By assumption, there is an $i$ such
  that $r\to x_i$ is not surjective or $r\to x^i$ is not
  injective. Since therefore $b$ or $a'$ is not an isomorphism,
  $(z_i)_{y_i}$ is a proper subquotient object of $x_i$. Since $\<r\>$
  factors through $[\z]_\y$ with $\z=x^i\times z_i$ and
  $\y=x^i\times y_i$ its image lies in $X$.

  Next we claim that
  \[\label{eq:HomX}
    \phi:\<T'\>_\KK\to\|Hom|_\cT(\1,X)
  \]
  is bijective. This will prove the theorem since then $\phi$ will map
  $\<T(x_1,\ldots,x_n)\>$ isomorphically onto $\bigoplus_i[x_i]^0$.

  Since injectivity follows from \eqref{eq:phiiso} it suffices to show
  surjectivity. Let $(z_i)_{y_i}$ be a subquotient object of $x_i$ for
  each $i$ and assume that least one is proper. Put
  $\z=\prod_iz_i$. Then we have a commutative diagram
  \[
    \cxymatrix{
      \<\sO(\z)\>_\KK\ar[r]^>>>{\phi_z}_>>>\sim\ar[d]_{\iota}&
      \|Hom|_\cT(\1,[\z])\ar[d]_{\iota'}\\
      \<\sO(\x)\>_\KK\ar[r]^>>>{\phi}_>>>\sim& \|Hom|_\cT(\1,[\x])}
  \]
  where both $\iota$ and $\iota'$ are componentwise induced by the
  diagrams
  \[
    \cxymatrix{&y_i\sur[ld]\inj[rd]\\\
      z_i&&x_i}
  \]
  Observe that $\iota$ has its image in $\<T'\>_\KK$. In fact, let
  $(z_i)_{y_i}$ be the proper subquotient object of $x_i$. Then
  diagram \eqref{eq:xizixi} with $r$ replaced by $\iota(r)$ shows that
  $\iota(r)\to x_i$ is not surjective or $\iota(r)\to x^i$ is not
  injective, hence $\iota(r)\in T'$. Since $\phi_z$ is an isomorphism
  and $X$ is generated by the images of $\iota'$ the claim follows.
\end{proof}

\begin{remark}
  For $n\le1$ the theorem says $\|Hom|_\cT(\1,[\*]^0)=\KK$ (which is
  trivial since $[\*]^0=[\*]=\1$) and $\|Hom|_\cT(\1,[x_1]^0)=0$ for
  $x_1\not\cong\*$. For $n=2$ we get
  $\|Hom|_\cT(\1,[x_1]^0\otimes[x_2]^0)=0$ if $x_1\not\cong x_2$ and
  $\|Hom|_\cT(\1,[x]^0\otimes[x]^0)=\KK[\|Aut|(x)]$ in accordance with
  \cref{lemma:defx1}\ref{it:defx1-3}.
\end{remark}

Now put $A(x):=\|Aut|_\cT(x)$ and $A:=\prod_iA(x_i)$. By abuse of
notation let $\chi_T=\chi_{T(x_1,\ldots,x_n)}$ be the character of the
$A$-module $\<T(x_1,\ldots,x_n)\>_\KK$. Since it is a permutation
character we get the explicit formula
\[\label{eq:Tchar}
  \chi_T(g)=|T(x_1,\ldots,x_n)^g|=|\{r\in T(x_1,\ldots,x_n)\mid
  gr=r\}|.
\]

\begin{corollary}
  For $n\ge0$ let $x_1,\ldots,x_n$ be objects of $\cA$ and
  $\chi_1\otimes\ldots\otimes\chi_n$ an irreducible character of
  $A$. Then
  \[\label{eq:multtensor1}
    \begin{split}
      \dim_\KK\|Hom|&(\1,[x_1]^0_{\chi_1}\otimes\ldots\otimes[x_n]^0_{\chi_n})
      \\&=\<\chi_T\mid\chi_1\otimes\ldots\otimes\chi_n\>_A=
      \sum_{r\in T/A}\<\chi_{\|triv|}\mid\chi_1\otimes\ldots\otimes\chi_n\>_{A_r}=\\
                    &=\frac1{|A|}\sum_{g\in A}|T^g|\,\chi_1(g_1)\ldots\chi_n(g_n)=
                      \frac1{|A|}\sum_{\ueber{r\in T,g\in A}{ gr=r}}\chi_1(g_1)\ldots\chi_n(g_n)
    \end{split}
  \]
  where $g=(g_1,\ldots,g_n)\in A$, $T=T(x_1,\ldots,x_n)$, $T/A$ is a
  set of orbit representatives.
\end{corollary}

\begin{proof}
  By \eqref{eq:decomp} we have
  \[
    \|Hom|_\cT(\1,\bigotimes\nolimits_i[x_i]^0)=
    \bigoplus_{\chi_1,\ldots,\chi_n}
    \|Hom|_\cT(\1,\bigotimes\nolimits_i[x_i]_{\chi_i}^0)\otimes
    \bigotimes_iV_{\chi_i}.
  \]
  Hence \eqref{eq:abc2} implies that
  $\dim_\KK\|Hom|_\cT(\1,\bigotimes\nolimits_i[x_i]_{\chi_i}^0)$ is
  the multiplicity of $\bigotimes_iV_{\chi_i}$ in $\<T\>_\KK$. The
  decomposition of $T$ into orbits gives
  $\chi_T=\sum_{r\in T/A}\|ind|_{A_r}^A\chi_{\|triv|}$ where $A_r$ is
  the stabilizer of $r$ in $A$. The second formula follows. The third
  formula follows from \eqref{eq:Tchar} and the fourth is just a
  reformulation of the third.
\end{proof}

The case $n=3$ corresponds to a multiplicity formula for the tensor
product of simple objects. In the next statement we let $nX$ stand for
$X^{\oplus n}$.

\begin{corollary}\label{cor:multtensor2}
  Let $x_1$ and $x_2$ be objects of $\cA$. Then
  \[\label{eq:multtensor2}
    [x_1]^0_{\chi_1}\otimes[x_2]^0_{\chi_2}\cong \bigoplus_{x,\chi}\
    \<\chi_{T(x_1,x_2,x)}\mid \chi_1\otimes\chi_2\otimes\chi^\vee\>\
    [x]^0_\chi
  \]
  where $x$ runs through all subquotients of $x_1\times x_2$ and
  $\chi$ through all characters of $A(x)$.
\end{corollary}

\begin{proof}
  Let
  $[x_1]^0_{\chi_1}\otimes[x_2]^0_{\chi_2}\cong \sum_{x,\chi}\
  m_{\chi_1,\chi_2}^\chi\ [x]^0_\chi$. Observe that
  $([x]^0_\chi)^\vee=[x]^0_{\chi^\vee}$. Now application of
  $X\mapsto\dim\|Hom|_\cT(\1,X\otimes[x]^0_{\chi^\vee})$ and
  \eqref{eq:multtensor1} yield
  $\<\chi_T\mid
  \chi_1\otimes\chi_2\otimes\chi^\vee\>=m_{\chi_1,\chi_2}^\chi$. Finally,
  in order for $T\ne\emptyset$ there must exist
  $r\subseteq x_1\times x_2\times x$ with $r\into x_1\times x_2$ and
  $r\auf x$ which implies that $x$ is a subquotient of
  $x_1\times x_2$.
\end{proof}
  
\begin{example}
  We verify formula \eqref{eq:multtensor2} in the case of
  $\cA=\mathsf{Set}^{\|op|}$. Let $x=(n):=\{1,\ldots,n\}$. The
  irreducible representations $V_\lambda$ of $A(x)=S_n$ are in
  $1-1$-correspondence with partitions $\lambda$ of $n$. In
  particular, for every $\lambda$ there is a simple object
  $\{\lambda\}:=[(n)]_\lambda^0$.

  To the irreducible character $\chi_\lambda$ is assigned the Schur
  function $s_\lambda$. Then the tensor product corresponds to the
  \emph{Kronecker product} of symmetric functions:
  \[
    \text{When }V_\lambda\otimes V_\mu=\bigoplus_\nu
    m^\nu_{\lambda\mu}V_\nu \quad\text{ then }\quad
    s_\lambda*s_\mu:=\sum_\nu m^\nu_{\lambda\mu}s_\nu.
  \]
  Likewise, the tensor product of simple objects of $\cT$ gives rise
  to a product, as well:
  \[\label{eq:abc1}
    \text{When }\{\lambda\}\otimes \{\mu\}=\bigoplus_\nu
    M^\nu_{\lambda\mu}\{\nu\}\quad \text{ then }\quad s_\lambda\star
    s_\mu:=\sum_\nu M^\nu_{\lambda\mu}s_\nu.
  \]
  This product was introduced by Murnaghan \cites{Murnaghan,Murnaghan2}
  as \emph{stable Kronecker product} in another form: For a partition
  $\lambda$ and $n>0$ let $\tilde\lambda:=(n-|\lambda|,\lambda)$. It
  is a partition of $n$ as soon as $n\ge |\lambda|+\lambda_1$. Then we
  have the stability formula
  \[
    m^{\tilde\tau}_{\tilde\lambda\tilde\mu}=M^\tau_{\lambda\mu}\quad\text{for
    }n>\!\!>0.
  \]
  That these two definitions coincide follows from
  \cite{Deligne}*{Prop.~6.4} (see also \cite{EntovaAizenbud}).

  Using the scalar product on symmetric functions for which the Schur
  functions form an orthonormal basis one defines the (double) skew
  Schur functions by
  \[
    s_{\lambda\setminus\mu\nu}:=\sum_\tau\<s_\lambda|s_\mu s_\nu
    s_\tau\>s_\tau.
  \]
  We will show in \cref{sec:AppB} that \cref{cor:multtensor2} is
  equivalent to the following result of Littlewood:

\begin{theorem}[Littlewood \cite{Littlewood}*{Thm.~IX}]\label{thm:Littlewood}
  Let $\lambda$ and $\mu$ be partitions. Then
  \[
    s_\lambda\star s_\mu=
    \sum_{\ueber{\alpha,\beta,\gamma}{|\alpha|=|\beta|}}(\,s_\alpha*s_\beta)\,
    s_{\lambda\backslash\alpha\gamma}\, s_{\mu\backslash\beta\gamma}.
  \]
\end{theorem}

\end{example}

Using formula \eqref{eq:multtensor2}, we are able to determine the
Grothendieck ring $K(\cT)$. Recall that $\hat\cA$ denotes the set of
isomorphism classes of objects of $\cA$. The representation ring of a
finite group $G$ is denoted by $\cR(G)$. It will be considered as a
subring of the ring $\KK[G]^G$ of $\KK$-valued class functions.

\begin{corollary}
  As an additive group
  \[
    K(\cT)=\bigoplus_{x\in\hat\cA}\cR(\|Aut|_\cT(x)).
  \]
  Moreover, the product of $f_1\in\cR(A(x_1))$ and $f_2\in\cR(A(x_2))$
  is given by $f_1\star f_2=\sum_xf_1\star _xf_2$ with $x\in\hat\cA$
  running over all subquotients of $x_1\times x_2$ and
  \[
    f_1\star _xf_2=\<\chi_{T(x_1,x_2,x)}\mid f_1\otimes
    f_2\>_{A(x_1)\times A(x_2)}\in \cR(A(x)).
  \]

\end{corollary}

\begin{proof}
  Let $a_i:=|A(x_i)|$. The simple object $[x]^0_\chi$ corresponds to
  the function $\chi\in\cR(A(x))$. Hence, \eqref{eq:multtensor2}
  implies
  \[
    \begin{split}
      (f_1\star _xf_2)(h)&=
                           \frac1{a_1a_2|A(x)|}
                           \sum_{\chi}\ \<\chi_{T(x_1,x_2,x)}\mid
                           \chi_1\otimes\chi_2\otimes\chi^\vee\>\chi(h)=\\
                         &=\frac1{a_1a_2}\sum_{g_1,g_2,g}
                           \chi_T(g_1,g_2,g)f_1(g_1)f_2(g_2)
                           \frac1{|A(x)|}\sum_\chi\chi(g^{-1})\chi(h)=\\
                         &=\frac1{a_1a_2}\sum_{g_1,g_2}
                           \chi_T(g_1,g_2,h)f_1(g_1)f_2(g_2)=\\
                         &=\<\chi_T(*,*,h)\mid f_1\otimes f_2\>_{A(x_1)\times A(x_2)}.\hspace{150pt}\qedhere
    \end{split}
  \]
\end{proof}

In general, the components $f_1\star _xf_2$ are difficult to compute,
except for the highest one corresponding to $x=x_1\times x_2$. Observe
that in that case $A(x_1)\times A(x_2)$ is a subgroup of $A(x)$.

\begin{proposition}
  Let $x=x_1\times x_2$ Then
  \[
    f_1\ast_xf_2=\|ind|_{A(x_1)\times A(x_2)}^{A(x_1\times
      x_2)}f_1\otimes f_2.
  \]
\end{proposition}

\begin{proof}
  Let $r\in T(x_1,x_2,x)$. Since $r\to x^3=x_1\times x_2=x$ is
  injective and $r\to x_3=x$ is surjective, $r$ realizes $x$ as a
  subquotient of itself. By \cite{TERC}*{Lemma\ 2.6}, both morphisms
  are isomorphisms. So $r$ is the graph $=\Gamma_h$ of an isomorphism
  $h:x_1\times x_1\overset\sim\to x$. Let $g_1\in A(x_1)$,
  $g_2\in A(x_2)$ and $g\in A(x)$. Then $\Gamma_h$ is fixed by
  $(g_1,g_2,g)$ if and only if $h^{-1}gh=g_1\times g_2$. This implies
  \[
    \chi_T(g_1,g_2,g)=|T^{(g_1,g_2,g)}|=|\{h\in A(x)\mid
    h^{-1}gh=g_1\times g_2\}|
  \]
  and therefore
  \[
    f_1\star _xf_2(g)= \hbox to
    40pt{$\frac1{|A(x_1)\times
        A(X_2)|}$\hss} \sum_{\ueber{h\in A(x)}{ h^{-1}gh\in A(x_1)\times
      A(x_2)}}\hspace{-20pt}(f_1\otimes f_2)(h^{-1}gh)= \|ind|_{A(x_1)\times
      A(x_2)}^{A(x_1\times x_2)}f_1\otimes f_2.
  \]
\end{proof}

\begin{remark}
  In the case $\cA=\mathsf{Set}^{\|op|}$ this means that the top
  degree component of the stable Kronecker product is the usual
  product of symmetric functions
\end{remark}

Let $\Kq(\cT)=K(\cT)$ but with multiplication
$f_1f_2=f_1\ast_{x_1\times x_2}f_2$ for $f_i\in\cR(A(x_i))$. Then
$\Kq(\cT)$ becomes a commutative ring.

\begin{definition}
  An \emph{valuation} of $\cA$ is a function $v:\hat\cA\to\RR_{\ge0}$
  with the following properties:
  \begin{enumerate}
  \item $v(x_1\times x_2)=v(x_1)+v(x_2)$
  \item If $z\prec x$ then $v(z)<v(x)$.
  \end{enumerate}
\end{definition}

Every $\cA$ which has a forgetful functor $F:\cA\to\mathsf{Sets}$
preserving products, injectivity, and surjectivity has a valuation
namely $v(x)=\log|F(x)|$.

A valuation induces a filtration on the group $K(\cT)$ by
\[
  K_{\le a}(\cT):=\bigoplus_{\ueber{x\in\hat\cA}{ v(x)\le a}}\cR(A(x))
\]
Since $f_1\star _xf_2\ne0$ implies $x\preccurlyeq x_1\times x_2$ and
therefore $v(x)\le v(x_1)+v(x_2)$ this is also a filtration of
rings. Let
$\|gr|_vK(\cT):=\bigoplus_{a\ge0}K_{\le a}(\cT)/K_{<a}(\cT)$ be the
associated graded ring. Then

\begin{corollary}
  Let $v$ be a valuation of $\cA$. Then $\|gr|_vK(\cT)=\Kq(\cT)$.
\end{corollary}

\section{Dimensions}

Since $\|End|_\cT\1=\KK$ one can define an \emph{internal trace} for
any endomorphism $\phi$ of an object $X$ of $\cT$ by
\[\label{eq:defdim}
  \cxymatrix{\1\ar[r]\ar@/_1.5pc/[rrrr]_{\|tr|_X\phi}&X\otimes
    X^\vee\ar[r]^{\phi\otimes1}&X\otimes X^\vee\ar[r]^\sim&
    X^\vee\otimes X\ar[r]&\1}.
\]
The \emph{internal dimension of an object $X$} is then defined as the
trace of the identity: $\dim X:=\|tr|_X\|id|_X$. More generally, for
an idempotent $p$ on $X$ we have $\dim p(X)=\|tr|_Xp$.

\begin{lemma}\label{cor:dim[x]*}
  Let $x$ be an object of $\cA$. Then
  \[
    \dim_\cT[x]=\delta(x),\quad\dim_\cT[x]^*=\omega_x,\quad
    \dim[x]^0=\sum_{z\in\qO(x)}\mu_\qO(x,z)\omega_z.
  \]
\end{lemma}

\begin{proof}
  For $X=[x]$ the diagram \eqref{eq:defdim} leads to
  $\dim_\cT[x]=\<\Delta x\doppelpfeil\*\>=\delta(x)$. The second
  formula is derived as follows:
  \[
    \dim[x]^*\overset{\eqref{eq:p*}} =\sum_{y\subseteq
      x}\mu_\sO(y,x)\|tr|p_y= \sum_{y\subseteq
      x}\mu_\sO(y,x)\|dim|_\cT[y]=\sum_{y\subseteq
      x}\mu_\sO(y,x)\delta(y) \overset{\eqref{eq:omegax}}=\omega_x.
  \]
  The last formula is obtained by Möbius inversion of
  \eqref{eq:subquotient3}.
\end{proof}

The simple objects $[x]^0_\chi$ appear in the isotypical decomposition
\eqref{eq:decomp} of $[x]^0$. To compute their dimension we need to
calculate its $A(x)$-character. But first we need to compute the
characters of $[x]$ and $[x]^*$.

Let $r\in\sO(x\times x)$ be a relation. Then we define its fixed point
object as the subobject
\[
  x^r:=r\cap\Delta x=r\Times_{x\times x} x\subseteq x.
\]
One checks easily that
\[\label{eq:rfix}
  \|tr|_x\<r\>=\dim_\cT[x^r]=\delta(x^r).
\]
If $r$ is the graph of an automorphism $g$ of $x$ then we write
$x^r=x^g$. It is easily seen that $x^g$ is also the maximal subobject
on which $g$ acts as identity.

When a finite group $A$ acts on an object $X$ of $\cT$ then one defines
the \emph{character} of this action as the $\KK$-valued class function
$\chi_X(g):=\|tr|_X(g)$.

\begin{proposition}\label{thm:charx*}
  Let $g\in A(x)$. Then
  \[\label{eq:char1}
    \chi_{[x]}(g)=\dim_\cT[x^g] \quad\text{and}\quad \chi_{[x]^*}(g)=
    \frac{\omega_x}{|A(x)|}\chi_{\|reg|}(g)=
    \begin{cases}
      \omega_x,&\text{if }g=\|id|_x;\\0,&\text{otherwise.}
    \end{cases}
  \]
\end{proposition}

\begin{proof}
  The first assertion follows from \eqref{eq:rfix}. For the second we
  use induction on $|\sO(x)|$. So we assume that the claim is correct
  for all proper subobjects of $x$. If $g=\|id|_x$ the assertion is
  just \cref{cor:dim[x]*}. Otherwise there is, according to the
  subobject decomposition \eqref{eq:SOD}, a $A$-equivariant
  decomposition
  \[
    [x]=
    \bigoplus_{y\subseteq x^g}[y]^*\oplus\bigoplus_{y\not\subseteq x^g}[y]^*=
    [x^g]\oplus[x]^*\oplus X
    \quad\text{with}\quad
    X=\bigoplus_{\ueber{y\subsetneq x}{ y\not\subseteq x^g}}[y]^*.
  \]
  Only the summands $[y]^*$ of $X$ with $gy=y$ enter the computation
  of $\|tr|_X(g)$. For these, we have $y^g=y\cap x^g\ne y$ and
  therefore $g|_y\ne\|id|_y$. The induction hypothesis implies
  $\|tr|_{[y]^*}(g)=0$. Hence
  \[
    \|tr|_{[x]^*}(g)=\|tr|_{[x]}(g)-\|tr|_{[x^g]}(g)-\|tr|_X(g)=
    \dim_\cT[x^g]-\dim_\cT[x^g]-0=0.\qedhere\hfill\qed
  \]
\end{proof}

Next, we compute the character of $[x]^0$ using equivariant Möbius
inversion on \eqref{eq:subquotient3}. The procedure is basically the
same as that of Assaf-Speyer \cite{AssafSpeyer}*{\S5} except that we
avoid the explicit use of simplicial complexes.

For $z\in\qO(x)$ and $n\ge0$ let $\|Ch|^n(x,z)$ be the $\KK$-vector
spaces spanned by all chains in $\qO(x)$
\[
  x=z_0>\ldots>z_n=z.
\]
These combine to a graded vector space
$\|Ch|^*(x,z):=\bigoplus_n\|Ch|^n(x,z)$. Moreover, let
\[
  \|Ch|^*(x):=\bigoplus_{z\in\qO(x)}\|Ch|^*(x,z)\otimes_\KK[z]^*.
\]
This is a $\ZZ$-graded object of $\cT$. Assume a group $A$ acts on
$x$.  Then $\|Ch|^*(x)$ is an object of the $A$-equivariant category
$\cT^A$. Let $K_A(\cT)$ be its Grothendieck group. The class of an
object $X$ of $\cT^A$ in $K_A(X)$ will be denoted by
$\{X\}_A$. Observe, that for any subgroup $A'\subseteq A$ there is a
natural induction homomorphism
$\|ind|_{A'}^A=\KK[A]\otimes_{\KK[A']}\bullet:K_{A'}(\cT)\to
K_A(\cT)$. The following is well-known when $\cT$ is the category of
vector spaces. We repeat the arguments for $\cT$.

\begin{lemma}
  Let $A$ act on the object $x$ of $\cA$. Then the following identity
  holds in $K_A(\cT)$:
  \[
    \{[x]^0\}_A=\chi(\|Ch|^*(x)):=\sum_n(-1)^n\{\|Ch|^n(x)\}_A.
  \]
\end{lemma}

\begin{proof}
  We argue by induction on $|q(x)|$ and assume the formula to be
  correct for all proper quotients $z$ of $x$. It follows from the
  subquotient decomposition \eqref{eq:subquotient3} of $[x]^*$ that
  \[
    \{[x]^0\}_A=\{[x]^*\}_A-\{\bigoplus_{x>z}[z]^0\}_A.
  \]
  Let $\qO'/A\subseteq\qO(x)\setminus\{x\}$ be a set of $A$-orbit
  representatives. Using the $A$-isomorphism
  \[\label{eq:indAAu}
    \bigoplus_{x>z}[z]^0= \bigoplus_{z\in\qO'/A}\|ind|_{A_z}^A[z]^0.
  \]
  the induction hypothesis yields
  \[
    \{\bigoplus_{x>z}[z]^0\}_A=
    \sum_{z\in\qO'/A}\|ind|_{A_z}^A\chi(\|Ch|^*(z))=
    \sum_{z\in\qO'/A}\chi(\|ind|_{A_z}^A\|Ch|^*(z))=
    \chi(\bigoplus_{x>z}\|Ch|^*(z)).
  \]
  By prepending $x$ to a chain we see that
  \[
    \|Ch|^*(x)=\|Ch|^0(x)\oplus\bigoplus_{x>z}\|Ch|^{*+1}(z)
  \]
  Because of $\|Ch|^0(x)=[x]^*$ this implies
  \[
    \{\bigoplus_{x>z}[z]^0\}_A=\{[x]^*\}_A-\chi(\|Ch|^*(x))
  \]
  and the assertion follows.
\end{proof}

We now derive two formulas for the $A(x)$-character of $[x]^0$.
Recall that $\|Ch|^*(x,z)$ can be turned into a complex with the
differential
\[
  \delta:\|Ch|^n(x,z)\to\|Ch|^{n-1}(x,z):
  (z_0>\ldots>z_n)\mapsto\sum_{i=1}^{n-1}(-1)^i
  (z_0>\ldots>\hat{z_i}>\ldots>z_n).
\]
The lattice $x/z:=\qO(x)_{\le z}$ is modular (see
\cref{thm:qmodular}). Therefore (see \cref{thm:hommod}) the homology
of $\|Ch|^*(x/z)$ is concentrated in degree $r:=\|rk|x/z$ and the
group $A(x)_z$ acts on $H_r(\|Ch|^*(x/z))$ with a character
$\St_{x/z}$.

Define moreover the \8inertia group\9 of $z\in\qO(x)$ as
\[
  A(x/z):=\{g\in A(x)_z\mid g|_z=\|id|_z\},
\]
let $\hat x$ be the socle (i.e., the meet of all atoms) of
$\qO(x)$, and let $(-1)^{x/z}:=(-1)^{\|rk|x/z}$. Then we have:

\begin{theorem}
  The character of $[x]^0$ is
  \[\label{eq:decomp2}
    \chi_{[x]^0}=\sum_{z\in\qO(x/\hat x)}
    \frac{(-1)^{x/z}\omega_z}{[A(x):A(x/z)]}\
    \|ind|_{A(x/z)}^{A(x)}\St_{x/z}.
  \]
\end{theorem}

\begin{proof}
  Write, as in \eqref{eq:indAAu},
  \[
    \chi(\|Ch|^*(x))= \sum_{z\in\qO/A(x)}\|ind|_{A(x)_z}^{A(x)}
    \chi(\|Ch|^*(x,z))\,\{[z]^*\}_{A(x)_z}.
  \]
  Thus
  \[
    \chi_{[x]^0}= \sum_{z\in\qO/A(x)}\|ind|_{A(x)_z}^{A(x)}
    (-1)^{x/z}\St_{x/z}\,\chi_{[z]^*}.
  \]
  Now \eqref{eq:char1} implies that
  \[
    \St_{x/z}\,\chi_{[z]^*}= \frac{|A(x/z)|}{|A(x)_z)|}\
    \omega_z\,\|ind|_{A(x/z)}^{A(x)_z}\St_{x/z}.
  \]
  The assertion follows with $q(x/\hat x)$ replaced by $q(x)$. Now
  recall that $H_*(x/z)=0$ unless $\mu_\qO(x,z)\ne0$ (see
  \cref{thm:hommod}) which, according to \cref{lemma:muvanish}, is
  equivalent to $z\le \|soc|\qO(x)=\hat x$.
\end{proof}

\begin{corollary}
  Let $\chi$ be an irreducible character of $A(x)$. Then
  \[
    \dim_\cT[x]_\chi^0=\sum_{z\in\qO(x/\hat x)}
    \frac{(-1)^{x/z}\omega_z}{[A(x):A(x/z)]} \
    \<\St_{x/z}\mid\chi\>_{A(x/z)}
  \]
\end{corollary}

\begin{example}
  Let $x$ be a simple object, i.e., $|\qO(x)|=2$. Since $A(x/x)=1$,
  $A(x/\*)=A(x)$, and $h_{x/z}$ is the trivial character we get
  \[
    \dim[x]^0_\chi=\frac{\deg\chi}{|A(x)|}\cdot
    \begin{cases}
      \omega_x-|A(x)|&\text{if }\chi=\chi_{\|triv|}\\
      \omega_x&\text{otherwise}\\
    \end{cases}
  \]
  For example if $p$ is a prime then
  $\dim[\ZZ_p]^0_{\chi_{\|triv|}}=\frac{t_p-1}{p-1}-1=\frac{t_p-p}{p-1}$
  where $t_p=\delta(\ZZ_p)$.
\end{example}

\begin{example}
  Assume, more generally, that $x$ has a unique minimal quotient $z$
  such that $\hat x=z$. Then
  \[
    \dim_\cT[x]^0_\chi=\frac{\deg\chi}{|A(x)|}\cdot
    \begin{cases}
      \omega_x-|A(x/z)|\,\omega_z,&\text{if }\chi|_{A(x/z)}\text{ is trivial,}\\
      \omega_x,&\text{otherwise}\\
    \end{cases}
  \]
  Let for example $x=S_3$, the symmetric group in the category of
  finite groups then $z=S_3/A_3=\ZZ_2$. We have $A(x)=S_3$ and
  $A(x/z)=A(x)$. Since $\omega_z=t_2-1$ and $\omega_{x\auf z}=t_3-3$
  (where $t_2=\delta(\ZZ_2)$ and $t_3=\delta(\ZZ_3)$) we get
  \[
    \dim_\cT[S_3]^0_\chi=
    \begin{cases}
      \frac1{6}(t_2-1)(t_3-9)&\chi=\chi_{\|triv|}\\
      \frac{1}{6}(t_2-1)(t_3-3)&\chi=\chi_{\|sign|}\\
      \frac{1}{3}(t_2-1)(t_3-3)&\deg\chi=2\\
    \end{cases}
  \]
\end{example}

Our second approach is to compute $\chi_{[x]^0}(g)$ using the Hopf
fixed point formula (see \cite{Sundaram}).

To state the result, let $\qO^g:=\qO(x)^g:=\{z\in\qO(x)\mid g(z)=z\}$
be the fixed point sublattice of $\qO(x)$. It is modular, as well. Let
$\mu_{\qO^g}$ be its Möbius function and let $\hat x_g$ be its socle.

If $z\in\qO^g$ then $g$ will act on $z$. Suppose $g$ acts on two
objects $z_1,z_2\in\qO^g$ trivially. Since $z_1\wedge z_2$ is a
subobject of $z_1\times z_2$ (namely the image of $x$ under the
diagonal morphism), the element $g$ will also act trivially on
$z_1\wedge z_2$. This shows that there is an object $x_g\in\qO^g$ such
that $g|_z=\|id|_z$ if and only if $z\ge x_g$.

\begin{theorem}
  The character of $[x]^0$ at $g\in A(x)$ is
  \[\label{eq:decomp3}
    \chi_{[x]^0}(g)= \sum_{\ueber{z\in\qO(x)^{g}}{ x_g\le z\le\hat x_g}}
    \mu_{\qO(x)^g}(x,z)\ \omega_z
  \]
  
\end{theorem}

\begin{proof}
  The action of $g$ on $\|Ch|^n(x)$ is a permutation representation
  inside $\cT$. Since the trace of $g$ on any orbit which is not a
  fixed points is zero we see that the trace of $g$ on $\|Ch|^n(x)$ is
  equal to that on $\|Ch|^n(x,z;g)$ where the latter is the subspace
  of $\|Ch|^n(x,z)$ spanned by the $g$-invariant chains. Thus we get
  \[
    \begin{split}
      \|tr|(g:\|Ch|^n(x))
      &=
        \|tr|\left(g:\bigoplus_{z\in\qO^g}
        \|Ch|^n(x,z;g)\otimes_\KK[z]^*\right)=\\
      &=\sum_{z\in\qO^g}\|tr|(g:\|Ch|^n(x,z;g))\
        \|tr|_{[z]^*}(g)=\\
      &\overset{\eqref{eq:char1}}=
        \sum_{\ueber{z\in\qO^g}{ x_g\le z}}
        \|tr|(g:\|Ch|^n(x,z;g))\ \omega_z.
    \end{split}
  \]
  Since $\|Ch|^n(x,z;g)$ is just the chain complex of $\qO^g$ and $g$
  acts trivially on it, its trace equals the dimension which is
  $\mu_{\qO^g}(x,z)$ (see \eqref{eq:Hopf}). Since this Möbius number
  is $0$ unless $z\le\hat x_g$, the assertion follows.
\end{proof}

\begin{remarks}
  \emph{i)} Equation \eqref{eq:decomp3} is more economical than
  \eqref{eq:decomp2} since $\hat x_g\le\hat x_{\|id|}=\hat x$ for all
  $g$. In fact, for any atom $z\in\qO^g$ let $z_0\in\qO$ be the join
  of all atoms $z'\in\qO$ with $z'\le z$. Since $z_0\in\qO^g$ and
  $z_0\le z$ it follows $z=z_0$. This proves that $\hat x_g$ is the
  join of atoms of $\qO$ and therefore $\hat x_g\le\hat x$.

  \emph{ii)} It may happen that $x_g\not\le\hat x_g$, see, e.g., the
  third conjugacy class of \cref{ex:FF3} below. In that case
  $\chi_{[x]^0}(g)=0$.
\end{remarks}

The second dimension formula reads:

\begin{corollary}\label{cor:second}
  \[
    \dim_\cT[x]_\chi^0=\<\chi_{[x]^0}\mid\chi\>=
    \frac1{|A(x)|}\sum_{g\in A(x)}\chi(g)\sum_{x_g\le z\le\hat x_g}
    \mu_{\qO(x)^g}(x,z)\ \omega_z.
  \]
\end{corollary}


\begin{example}\label{ex:FF3}
  Let $x=\FF_2^3$ in the category of $\FF_2$-vector spaces. Then
  $A(x)=GL(3,\FF_2)$ is the simple group of order $168$. Its six
  conjugacy classes are represented by the following matrices. Below
  we indicated the lattice $\qO(x)^g$ (or rather the kernels with
  $x_g$ in bold and $\hat x_g$ underlined) and
  $\chi_{[x]^0}(g)=\sum_{z\le x_g}\mu_{\qO(x)^g}(x,z)\omega_z$.
  
  $\begin{array}{cccccc}
     \begin{pmatrix}
       1&&\\
        &1&\\
        &&1
     \end{pmatrix}
        &
          \begin{pmatrix}
            1&1&\\
             &1&\\
             &&1
          \end{pmatrix}
        &
          \begin{pmatrix}
            1&1&\\
             &1&1\\
             &&1
          \end{pmatrix}
        &
          \begin{pmatrix}
            &&1\\
            1&&1\\
            &1&
          \end{pmatrix}
          &
            \begin{pmatrix}
              &&1\\
              1&&\\
              &1&1
            \end{pmatrix}
          &
            \begin{pmatrix}
              &1&\\
              1&1&\\
              &&1
            \end{pmatrix}
     \\
     \\
     \qO(\FF_2^3)
        &
          \Hide{\kxymatrix{
        &e_1e_2e_3\ar@{-}[dl]\ar@{-}[dr]\\
     e_1e_2&\underline{e_1e_3}
             \ar@{}[u]^(0){}="a"^(0.9){}="b" \ar@{-} "a";"b"&e_1e_{23}\\
     e_{13}\ar@{-}[ur]&\mathbf{e_1}\ar@{-}[ul]
                        \ar@{}[u]^(0){}="a"^(0.8){}="b" \ar@{-} "a";"b"
                        \ar@{-}[ur]&e_3\ar@{-}[ul]\\
        &0\ar@{-}[ul]\ar@{-}[ur]
          \ar@{}[u]^(0.1){}="a"^(0.9){}="b" \ar@{-} "a";"b"
          }}
        &
          \kxymatrix{
          e_1e_2e_3\\
     \mathbf{e_1e_2}\ar@{}[u]^(0){}="a"^(0.9){}="b" \ar@{-} "a";"b"\\
     \underline{e_1}\ar@{}[u]^(0.1){}="a"^(0.9){}="b" \ar@{-} "a";"b"\\
     0\ar@{}[u]^(0.1){}="a"^(0.8){}="b" \ar@{-} "a";"b"\\
     }
        &
          \kxymatrix{
          \underline{\mathbf{e_1e_2e_3}}\\
     \\
     0\ar@{-}[uu]\\
     }
        &
          \kxymatrix{
          \underline{\mathbf{e_1e_2e_3}}\\
     \\
     0\ar@{-}[uu]\\
     }
        &
          \Hide{\kxymatrix{
        &\hide{\underline{e_1e_2e_3}}\ar@{-}[ddr]\\
     \mathbf{e_1e_2}\ar@{}[ur]^(.05){}="a"^(.7){}="b" \ar@{-} "a";"b"\\
        &&e_3\\
        &0\ar@{-}[uul]\ar@{-}[ur]\\
     }}
     \\\\
     t^3{-}14t^2{+}49t{-}44&-(t-1)(t-4)&0&-1&-1&-(t-2)\\\\
   \end{array}$
    
   Here $e_{ij}:=e_i+e_j$ and $e_ie_j=\<e_i,e_j\>$.  This implies the
   following dimension formulas where $\chi_d$ signifies an
   irreducible character of degree $d$.
   \[
     \begin{array}{ll}
       \dim[\FF_2^3]^0_{\chi_1}\ =&\sfrac1{168}(t - 1)(t - 2)(t - 2^5)\\
       \dim[\FF_2^3]^0_{\chi_3}\ =&\sfrac3{168}(t - 1)(t - 2)(t - 2^2)\\
       \dim[\FF_2^3]^0_{\chi_3'}\ =&\sfrac3{168}(t - 1)(t - 2)(t - 2^2)\\
       \dim[\FF_2^3]^0_{\chi_6}\ =&\sfrac6{168}(t - 1)(t - 2^2)(t - 2^4)\\
       \dim[\FF_2^3]^0_{\chi_7}\ =&\sfrac7{168}(t - 1)(t - 2)(t - 2^3)\\
       \dim[\FF_2^3]^0_{\chi_8}\ =&\sfrac8{168}(t - 2)(t - 2^2)(t - 2^3)
     \end{array}
   \]
 \end{example}

  \begin{example}
    Let $\cA=\FS^{\|op|}$. We adopt the notation of \cite{Macdonald},
    i.e. let $\chi^\lambda$ be the irreducible character of $S_n$
    corresponding to a partition $\lambda$ of $n$. For another
    partition $\mu$ of $n$ let $\chi^\lambda_\mu$ be the value of
    $\chi^\lambda$ in a permutation $g_\mu$ of cycle type $\mu$. Let
    $m_i(\mu)=\mu_i'-\mu_{i+1}'$ be the number of parts equal to $i$,
    so $m_1(\mu)$ is also the number of fixed points of $g_\mu$. The
    order of the centralizer of $g_\mu$ is
    $z_\mu:=\prod_i(i^{m_i(\mu)}m_i(\mu)!)$.
  
    Let $x=(n):=\{1,\ldots,n\}$ be an object of $\FS^{\|op|}$. Then
    $x_{g_\mu}=(n)^{g_\mu}\cong(m)$ and
    $\qO(x)^{g_\mu}=\mathfrak P((n))^{g_\mu}\cong\mathfrak P((b))$ where
    $m=m_1(\mu)$ and $b=\ell(\mu)=\mu_1'$. Then
    \[
      \chi_{[x]^0}(g_\mu)=\sum_{Z\subseteq[n]^{g_\mu}}\mu_{\mathfrak
        P^{g_\mu}}(Z,[n])\omega_Z= \sum_{k=0}^m(-1)^{b-k}\begin{pmatrix}
        m\\k\end{pmatrix}(t)_k=(-1)^{\ell(\mu)}\cC_{m_1(\mu)}
    \]
    where $(t)_k:=t(t-1)\ldots(t-k+1)$ and
    \[
      \cC_m(t):=\sum_{k=0}^m(-1)^k
      \begin{pmatrix} m\\k\end{pmatrix}(t)_k
    \]
    is, up to a sign, the $m$-th Charlier polynomial for the parameter
    value $a=1$. Combining this with a calculation of Deligne one
    obtains the following identity:

  \begin{corollary}\label{cor:identity}
    Let $f(x):=(t-x_1)\ldots(t-x_n)$ and $\delta:=(n-1,\ldots,0)$,
    Then
    \[\label{eq:identity}
      \dim[n]^0_{\chi^\lambda}=\sum_\mu(-1)^{\ell(\mu)}
      z_\mu^{-1}\chi^\lambda_\mu\cC_{m_1(\mu)}=
      \frac{\deg\chi^\lambda}{n!}f(\lambda+\delta).
    \]
  \end{corollary}
  
  \begin{proof}
    The first equality follows from
    $\dim[n]^0_{\chi^\lambda}=\<\chi_{[n]^0}\mid\chi^\lambda\>$. The
    right hand term is the dimension of $[x]^0_{\chi^\lambda}$ as
    computed by Deligne \cite{Deligne}*{7.4}.
  \end{proof}
  
  \begin{remark}
    In \cref{sec:AppB} we give a non-category theoretical proof of the
    right hand equality thereby obtaining an independent proof of Deligne's
    formula.
  \end{remark}

\end{example}

We conclude this paper with an interesting observation namely that in
each example the dimension of a simple object has nice
factorization. We show now that this is no coincidence. Let
$\KK_0\subseteq\KK$ be the $\overline\QQ$-algebra generated by all
values $\delta(e)$ with $e$ surjective and indecomposable. Then we have:

\begin{theorem}
  Let $S$ be a simple object of $\cA$. Then there are indecomposable
  surjective morphisms $e_1,\ldots,e_n$ (not necessarily distinct)
  such that $\dim_\cT S$ divides $\omega_{e_1}\ldots\omega_{e_n}$
  within $\KK_0$.
\end{theorem}

\begin{proof}
  Let $\KK_1$ be $\KK_0$ with all inverses $\omega_e^{-1}$, $e$
  surjective, adjoined.  Let $\KK'$ be another field and let
  $\phi:\KK_1\to\KK'$ be a ring homomorphism. Then
  $\delta'=\phi\circ\delta$ is by construction a non-degenerate degree
  function. Hence $\cT':=\cT(\cA,\delta')$ is a semisimple tensor
  category. Let now $S=[x]^0_\chi$ be the simple object of $\cT$. Then
  $S'=[x]^{0\prime}_\chi$ the corresponding simple object of
  $\cT'$. It follows from any of the explicit formulas above that
  first of all $d:=\dim_\cT S$ is an element of $\KK_0$ and that
  $\dim_{\cT'}S'=\phi(d)$. According to
  \cite{Deligne}*{Prop.~5.7(iii)}, the dimension of any simple object
  of $\cT'$ is nonzero. This implies $\phi(d)\ne0$ for all possible
  $\phi$ and therefore that $d$ is invertible in $\KK_1$. Thus there
  is $a\in\KK_0$ and $\omega=\omega_{e_1}\ldots\omega_{e_n}\in\KK_0$
  with $\frac a\omega\cdot d=1$ which shows $ad=\omega$ as asserted.
\end{proof}

\begin{examples}
  \emph{i)} For $\cA=\mathsf{Set}^{\|op|}$ and $t=\delta(\{0\})$
  transcendental over $\QQ$ we have $\KK_0=\overline\QQ[t]$ and
  $\omega=t-a$, $a\in\NN$ for $e$ indecomposable. Thus $\dim_\cT S$ is
  a product of linear factors $t-a$, $a\in\NN$. This follows of course
  also from Deligne's formula for $S$.

  \emph{ii)} For $\cA=\mathsf{Vect}(\FF_q)$ and $t=\delta(\FF_q)$
  transcendental over $\overline\QQ$ we have $\KK_0=\overline\QQ[t]$
  and $\omega=t-q^a$, $a\in\NN$ for every indecomposable
  $\omega_e$. Thus $\dim_\cT S$ is a product of linear factors
  $t-q^a$, $a\in\NN$.
\end{examples}
    
\begin{appendix}

\section{The Möbius function of the lattice of quotient objects}

In this section, we review some facts for modular lattices and their
relation to Mal'cev categories. The reason for our interest is the
following theorem of Pedicchio:

\begin{theorem}[\cite{Pedicchio}*{Cor.~2.3}]\label{thm:qmodular}
  Let $x$ be an object of a regular exact Mal'cev category. Then its
  set of quotients $\qO(x)$ is a modular lattice.
\end{theorem}

Recall that a lattice $\Lambda$ is \emph{modular} if is obeys the
modular identity
\[
  x\wedge(y\vee z)=(x\wedge y)\vee z\text{ whenever }x\le z.
\]
It is clear that with $\Lambda$ also every interval
\[
  x/z:=\{y\in\Lambda\mid x\le y\le z\}
\]
is modular, as well. The same holds for every sublattice.

In a modular lattice, all maximal chains
\[
  \hat0=z_0<z_1<\ldots<z_n=\hat1
\]
have the same length. This number $n$ is the \emph{rank}
$\|rk|\Lambda$ of $\Lambda$.

We are mostly interested in the Möbius function of a modular
lattice. Assume therefore that $\Lambda$ is finite with minimum
$\hat0$ and maximum $\hat1$. Recall that an \emph{atom} is an element
$z\in\Lambda$ with $|\hat0/z|=2$. The \emph{socle} $\|soc|(\Lambda)$
of $\Lambda$ is the join of its atoms. Recall also that a lattice is
\emph{complemented} if for every $x\in\Lambda$ there is $x'\in\Lambda$
with $x\wedge x'=\hat0$ and $x\vee x'=\hat1$. It is known that if a
modular lattice is complemented then all of its intervals are
complemented, as well.

\begin{lemma}\label{lemma:muvanish}
  Let $\Lambda$ be a finite modular lattice. Then for $z\in\Lambda$
  the following are equivalent:
  \begin{enumerate}[noitemsep]
  \item\label{it:4} $\mu(\hat0,z)\ne0$.
  \item\label{it:1} $z$ is the join of atoms.
  \item\label{it:2} $[\hat0,z]$ is complemented.
  \item\label{it:5} $z\le\|soc|(\Lambda)$.
  \end{enumerate}
\end{lemma}

\begin{proof}
  \ref{it:4}$\Rightarrow$\ref{it:1} See \cite{Stanley}*{3.9.5}.
  \ref{it:1}$\Rightarrow$\ref{it:2} See \cite{Birkhoff}*{IV \S5,
    Thm. 6}.
  \ref{it:2}$\Rightarrow$\ref{it:4} This is proved in the same way as
  \cite{Stanley}*{Prop.~3.10.1} together with the remark that in a
  complemented lattice the sum in \emph{loc.cit. (3.33)} is not empty.
%
%
%
  \ref{it:1}$\Rightarrow$\ref{it:5} Obvious.
  \ref{it:5}$\Rightarrow$\ref{it:2} The interval
  $[\hat0,\|soc|(\Lambda)]$ is complemented by the implication
  \ref{it:1}$\Rightarrow$\ref{it:2}. Thus also $[\hat0,z]$ is
  complemented.
\end{proof}

Even though it is not going to be used in this paper it may of
interest to know that finite complemented modular lattice have been
classified by Birkhoff \cite{Birkhoff}:

\begin{theorem}
  Every finite complemented modular lattice is product of
  indecomposable finite complemented modular lattices with factors
  which are unique up to order. Every indecomposables is isomorphic to
  one of the following lattices:
  \begin{itemize}[noitemsep]
  \item[$\|rk|=1$]The two element Boolean lattice
    $B_1:=\{\hat0,\hat1\}$.
  \item[$\|rk|=2$] The lattices
    $M_{q+1}:=\{\hat0<a_0,\ldots,a_q<\hat1\}$, $q\ge2$.
  \item[$\|rk|=3$] The subspace lattice of a finite projective plane.
  \item[$\|rk|\ge4$] The subspace lattice $L_n(q)$ of $\FF_q^n$.
  \end{itemize}
\end{theorem}

Using \cite{Stanley}*{(3.34) and Prop.~3.8.2.} one obtains:

\begin{corollary}\label{cor:mobius}
  Let $\Lambda$ be a finite complemented modular lattice. If $\Lambda$
  is indecomposable then
  \[
    \mu_\Lambda(\hat0,\hat1)=(-1)^rq^{\binom r2}
  \]
  where $r=\|rk|\Lambda$ and $q$ is as above. When
  $\Lambda=\prod_i\Lambda_i$ then
  \[
    \mu_\Lambda(\hat0,\hat1)=\prod_i\mu_{\Lambda_i}(\hat0,\hat1).
  \]
\end{corollary}

\begin{remark}
  The subspace lattice of $\FF_q^3$ is a special case of the
  $\|rk|=3$-case but there exists a plethora of others. Therefore it
  is a nice surprise (by Jónsson \cite{Jonsson}) that the latter do
  not occur as quotient lattice $\qO(x)$, at least when $x$ is a
  model for a Mal'cev algebraic theory. Likewise, in the same
  context, the lattices $M_m$ appear only when $q$ is a prime power
  (use \cite{MMT}*{Lemma~4.153 and Thm.~4.155}). Thus for algebraic
  theories only subspace lattices of vector spaces $\FF_q^n$, $n\ge1$,
  occur. Presumably this holds for all of our categories $\cA$.
\end{remark}

It is well known that Möbius numbers are categorified by the homology
of the complex of chains. Let $\Lambda$ be a finite lattice For
$n\ge0$ let $\|Ch|^n(\Lambda)$ be the $\KK$-vector space spanned by
all chains
\[
  \hat0=a_0<\ldots<a_n=\hat1.
\]
These spaces are connected by the differential
\[
  \delta:\|Ch|^n(\Lambda)\to\|Ch|^{n-1}(\Lambda):
  (a_0<\ldots<a_n)\mapsto\sum_{i=1}^{n-1}(-1)^i
  (a_0<\ldots<\hat{a_i}<\ldots<z_n).
\]

\begin{theorem}\label{thm:hommod}
  Let $\Lambda$ be a finite lattice. Then $(\|Ch|^*(\Lambda),\delta)$
  is a chain complex. If $\Lambda$ is modular of rank $r$ then its
  homology is concentrated in degree $r$ and the homology group
  $H(\Lambda)$ in degree $r$ has dimension
  $(-1)^r\mu_\Lambda(\hat0,\hat1)$.
\end{theorem}

\begin{proof}
  See e.g. \cite{Baclawski}*{Prop.~4.2 and Prop.~3.5}.
\end{proof}

If a group $A$ acts on $\Lambda$ then the character of $H(\Lambda)$
will be denoted by $h_\Lambda$. According to the Hopf trace formula we
have
\[\label{eq:Hopf}
  h_\Lambda(g)=(-1)^{\|rk|\Lambda}\mu_{\Lambda^g}(\hat0,\hat1).
\]
See e.g. \cite{Sundaram}*{(1.3)}. Here $\Lambda^g\subseteq\Lambda$ is
the fixed point sublattice of $g\in A$. Since $\Lambda^g$ is also
modular, this implies $h_\Lambda(g)=0$ unless the action of $g$ is
\8semisimple\9, i.e., every $g$-stable $a\in\Lambda$ has a $g$-stable
complement. In that case, $\Lambda^g$ is complemented and the Möbius
number can be computed using \cref{cor:mobius}. Note that when one
does these computations for $\Lambda=L_n(q)$ and $g\in GL(n,\FF_q)$
one obtains simply the Steinberg character at $g$ (see
e.g. \cite{CLT}).

\section{Symmetric function identities}\label{sec:AppB}

In this appendic we provide proofs for two identities in the main
text.

\subsection{Proof of \cref{thm:Littlewood}}\

Let $X_1,X_2,X_3$ be finite sets with $a_i:=|X_i|$ elements. Thus
$A=S_{a_1}\times S_{a_2}\times S_{a_3}$. An element of
$T(X_1,X_2,X_3)$ is a finite set $R$ together with injective maps
$X_i\into R$ such that $R$ is the union of any two images. Thus, the
set $T/A$ is parameterized by the four numbers
\[
  n_0:=|X_1\cap X_2\cap X_3|\text{ and }n_i:=|(X_j\cap X_k)\setminus
  X_i|
\]
for all permutations $(i,j,k)$ of $(1,2,3)$. These numbers are subject
to the constraints
\[
  a_i=n_0+n_j+n_k.
\]
The stabilizer of $R$ in $A$ is
$A_R=S_{n_0}\times S_{n_1}\times S_{n_2}\times S_{n_3}$. Its embedding
into $A$ is the composition of two homomorphisms, namely
\[
  S_{n_0}\times S_{n_1}\times S_{n_2}\times S_{n_3}\into
  \begin{matrix}
    S_{n_0}\times S_{n_0}\times S_{n_0}\times\\
    S_{n_3}\times S_{n_1}\times S_{n_2}\times\\
    S_{n_2}\times S_{n_3}\times S_{n_1}\phantom{\times}
  \end{matrix}
  \into S_{a_1}\times S_{a_2}\times S_{a_3}
\]
where the left arrow is four times a diagonal embedding and the right
embedding is columnwise.

We calculate the character of
$\<A/A_R\>_\KK=\|ind|_{A_R}^A\chi_{\|triv|}$ in terms of Schur
functions in two steps. First recall the formulas
\[
  \|ind|_{S_n}^{S_n\times S_n}\chi_{\|triv|}=\sum_{|\lambda|=n}
  s_\lambda\otimes s_\lambda
\]
and
\[
  \|ind|_{S_n}^{S_n\times S_n\times
    S_n}\chi_{\|triv|}=\sum_{|\alpha|=|\beta|=n} s_\alpha\otimes
  s_\beta\otimes (s_\alpha*s_\beta).
\]
Recall also that
\[
  \|ind|_{S_m\times S_n}^{S_{m+n}}s_\lambda\otimes s_\mu=s_\lambda
  s_\mu.
\]
Then
\[\label{eq:abc3}
  \|ind|_{A_R}^A\chi_{\|triv|}=
  \sum_{\alpha,\beta,\gamma_1,\gamma_2,\gamma_3} s_\alpha s_{\gamma_3}
  s_{\gamma_2}\otimes s_\beta s_{\gamma_1} s_{\gamma_3}\otimes
  (s_\alpha* s_\beta) s_{\gamma_2} s_{\gamma_1}
\]
where the summation variables are partitions subject to the
constraints
\[
  |\alpha|=|\beta|=n_0\text{ and }|\gamma_i|=n_i.
\]
For three partitions $\lambda$, $\mu$, $\nu$ with $|\lambda|=a_1$,
$|\mu|=a_2$, and $|\nu|=a_3$ we deduce from equations \eqref{eq:abc1},
\eqref{eq:abc2}, and \eqref{eq:abc3}
\begin{align}\nonumber
  \<s_{\lambda}\star s_{\mu}\mid
  s_{\nu}\>&{}=M_{\lambda\mu}^{\nu}=\sum_{R\in T/A}\<\|ind|_{A_R}^A\chi_{\|triv|}\mid
             s_{\lambda}\otimes s_{\mu}\otimes s_{\nu}\>=\\\label{eq:abc4}
           &{}=\sum_{\alpha,\beta,\gamma_1,\gamma_2,\gamma_3}
             \<s_\alpha s_{\gamma_3} s_{\gamma_2}\mid s_{\lambda}\>\ 
             \<s_\beta s_{\gamma_1} s_{\gamma_3}\mid s_{\mu}\>\ 
             \<(s_\alpha* s_\beta) s_{\gamma_2} s_{\gamma_1}\mid s_{\nu}\>
\end{align}
Only the constraint $|\alpha|=|\beta|$ needs to be kept since all
other possibilities contribute summands equal to $0$. After
rearranging the terms of \eqref{eq:abc4} we get the assertion:
\[
  \begin{split}
  \<s_{\lambda}\star s_{\mu}\mid
  s_{\nu}\>&=\sum_{\ueber{\alpha,\beta,\gamma_3}{|\alpha|=|\beta|}}
              \Big\<(s_\alpha* s_\beta)
              \ \underbrace{\textstyle\sum\limits_{\gamma_2}
              \<s_\alpha s_{\gamma_3} s_{\gamma_2}\mid s_{\lambda}\> s_{\gamma_2}
              }_{s_{\lambda\setminus\alpha\gamma_3}}
              \ \underbrace{\textstyle\sum\limits_{\gamma_1}
              \<s_\beta s_{\gamma_1} s_{\gamma_3}\mid s_{\mu}\> s_{\gamma_1}
              }_{s_{\mu\setminus\beta\gamma_3}}
              \ \Big|\  s_{\nu}\Big\>\\
            &=\<\sum_{\ueber{\alpha,\beta,\gamma_3}{|\alpha|=|\beta|}}
              (s_\alpha*s_\beta)s_{\lambda\setminus\alpha\gamma_3}
              s_{\mu\setminus\beta\gamma_3}\mid s_{\nu}\>\hspace{150pt}\qed
  \end{split}
\]

\subsection{Proof of \cref{cor:identity}}\

Multiplying all terms of \eqref{eq:identity} by $\chi_\tau^\lambda$
and summing over $\lambda$ one obtains the equivalent formula
\[\label{cor:identity2}
  \chi_{[n]^0}(g_\mu)=(-1)^{\ell(\mu)}\cC_{m_1(\mu)}=
  \sum_\lambda\sfrac{\deg\chi^\lambda}{n!}f(\lambda+\delta)\chi^\lambda_\mu.
\]
Consider the differential operators $\tilde D:=f(x_i\partial_{x_i})$
and $D:=a_\delta^{-1}\tilde D a_\delta$ (for unexplained notation see
\cite{Macdonald}*{pp. 467--468}). Then
$\tilde D(a_{\lambda+\delta})=f(\lambda+\delta)a_{\lambda+\delta}$ and
therefore $D(s_\lambda)=f(\lambda+\delta)s_\lambda$. From
$p_{1^n}=\sum_\lambda\deg\chi^\lambda s_\lambda$ and
$\<s_\lambda\mid p_\mu\>=\chi^\lambda_\mu$ one gets that the right
hand side of \eqref{cor:identity2} equals
$\frac1{n!}\<D(p_{1^n})\mid p_\mu\>$.

Now observe that, by Harish-Chandra's isomorphism, $D$ is the radial
part of a biinvariant differential operator $C$ on $M(n,\CC)$. Other
such operators are the Capelli operators
\[
  C_k:=\sum_{\ueber{I,J\subseteq[n]}{|I|=|J|=k}}\det A_{IJ}\det\partial
  A_{IJ},\quad k=0,\ldots,n.
\]
The radial part of $C_k$ is the polynomial
$P_{1^k}(x):=s_{1^k}^*(x-\delta)$ where $s_\lambda^*$ is a shifted
Schur function (see \cite{OO}*{Cor.~6.6}). From the relation
$f(x)=\sum_{k=0}^n(-1)^k(t)_{n-k}P_{1^k}$ (follows from
\cite{OO}*{(12.4)} by multiplying by $(u\downharpoonright r)$ and
setting $u=-t+n-1$) we get $C=\sum_{k=0}^n(-1)^k(t)_{n-k}C_k$.

The polynomial $p_{1^n}=e_1^n=(x_1+\ldots+x_n)^n$ is the radial part
of $h(A)=(\|tr|A)^n$ and we have $\det\partial A_{IJ}(h)=0$ for
$I\ne J$ and $\det\partial A_{II}(h)=(n)_{|I|}(\|tr|A)^{n-|I|}$. It
follows that the radial part of $C_k((\|tr|A)^n)$ is
$(n)_ke_ke_1^{n-k}$ and therefore
\[
  C(p_{1^n})=\sum_{k=0}^n(-1)^{n-k}\frac{n!}{k!}(t)_ke_{n-k}e_1^k.
\]
From $e_{n-k}=\sum_{|\tau|=n-k}(-1)^{n-k-\ell(\tau)}z_\tau^{-1}p_\tau$
(\cite{Macdonald}*{Ch.~I~(2.14$'$)}), $p_\tau e_1^k=p_{\tau1^k}$, and
$\<p_\lambda\mid p_\mu\>=z_\mu\delta_{\lambda\mu}$
(\cite{Macdonald}*{Ch.~I~(4.7)}) we get
\[
  \<\frac1{n!}D(p_{1^n})\mid p_\mu)= \sum_{k:\
    \tau1^k=\mu}(-1)^{\ell(\mu)-k}\frac{z_\mu}{z_\tau}\frac1{k!}(t)_k=
  (-1)^{\ell(\mu)}\sum_{k=0}^{m_1(\mu)}(-1)^k{\binom{m_1(\mu)}k}(t)_k
\]
as claimed.\qed

\end{appendix}

\end{document}